\documentclass{amsart}
\usepackage{geometry}
\usepackage{graphicx}
\usepackage{amssymb}
\usepackage[usenames,dvipsnames]{color}

\usepackage{hyperref}
\usepackage{mydef,boldfonts,xspace}

\makeatletter

\newcommand{\dxx}{\partial_{xx} }
\newcommand{\dyy}{\partial_{yy} }
\newcommand{\dzz}{\partial_{zz} }
\newcommand{\dyx}{\partial_{yx} }
\newcommand{\dxy}{\partial_{xy} }

\makeatother
\begin{document}

\title[Direction Splitting]{Convergence Analysis of a Class of Massively
  Parallel Direction Splitting Algorithms for the Navier-Stokes
  Equations}

\author[J.-L.~Guermond]{Jean-Luc Guermond$^1$} \address{$^1$Department
  of Mathematics, Texas A\&M University 3368 TAMU, College Station, TX
  77843-3368, USA.  On leave from CNRS, France.}
\email{guermond@math.tamu.edu}

\author[P.D.~Minev ]{Peter D. Minev $^2$} \address{$^2$Department of
  Mathematical and Statistical Sciences, University of Alberta,
  Edmonton, Alberta Canada T6G 2G1.}  \email{minev@ualberta.ca}

\author[A.J.~Salgado]{Abner J. Salgado$^3$} \address{$^3$Department of
  Mathematics, University of Maryland, College Park, MD 20742, USA.}
\email{abnersg@math.umd.edu}

\thanks{This material is based upon work supported by the National
  Science Foundation grants DMS-0713829 and a Discovery grant of the
  National Science and Engineering Research Council of Canada. This
  publication is also partially based on work supported by Award
  No. KUS-C1-016-04, made by King Abdullah University of Science and
  Technology (KAUST). The work of P.~Minev is also supported by
  fellowships from the Institute of Applied Mathematics and
  Computational Science and the Institute of Scientific Computing at
  Texas A\&M University. A.J.~Salgado is also supported by NSF grants
  CBET-0754983 and DMS-0807811.}

% \date{Draft Version: \today}

\keywords{Navier-Stokes, Fractional Time-Stepping, Direction Splitting}

\subjclass[2000]{
65N12,     % Stability and convergence of numerical methods
65N15,     % Error bounds
35Q30.    % PDE, Navier Stokes
}

\begin{abstract}
  We provide a convergence analysis for a new fractional
  time-stepping technique for the incompressible Navier-Stokes
  equations based on direction splitting. This new technique is of
  linear complexity, unconditionally stable and convergent, and
  suitable for massive parallelization.
\end{abstract}

\maketitle

\section{Introduction}
\label{sec:Intro}
This work is concerned with the analysis of a new class of
approximation techniques for the solution of the time-dependent
incompressible Navier-Stokes equations based on direction
splitting. This new technique requires, independently of the space
dimension, only the solution of a sequence of one-dimensional
problems, thus having linear complexity. The main claims of this paper
are that this technique is unconditionally stable and
superlinearly convergent with respect to the time discretization
parameter and is suitable for massive parallelization.

We consider the Stokes equations written in terms of velocity $\ue$
and pressure $\pe$ on a finite time interval $[0,T]$ and in a cubic
domain $\Omega = (0,1)^d$ with $d=2$ or $3$:
\begin{equation}
  \begin{cases}
    \ue_t - \LAP \ue + \GRAD \pe = f, & \text{in } \Omega\times(0,T], \\
    \DIV \ue = 0, & \text{in } \Omega\times[0,T], \\
    \ue|_{\partial\Omega} = 0, & \text{in } (0,T], \\
    \ue|_{t=0} = \ue_0, & \text{in } \Omega.
  \end{cases}
\label{eq:Stokes}
\end{equation}
where $f$ is a smooth source term and $\ue_0$ is a solenoidal initial
velocity field with zero normal trace. The nonlinear term in the
momentum equation of the Navier-Stokes equations is not accounted for
since it does not interfere with the incompressibility constraint.
The fluid density is assumed to be constant and has been put into the
normalization constants.

Once time is discretized, \eqref{eq:Stokes} reduces to a generalized
Stokes system at each time step.  Solving this coupled system often
proves computer intensive and is not easy to solve efficiently in
parallel due to the saddle point structure induced by the
incompressibility constraint.  Alternative more efficient approaches
consist of uncoupling the velocity and the pressure using so-called
projection algorithms.

Projection algorithms date back to the late 1960s and stem from the
seminal works of Chorin \cite{Chor68} and Temam \cite{tema69}. These
methods and various improvements thereof are still, to the best of our
knowledge, the methods of choice in the CFD community.  Although in
the 1980s and 1990s these techniques underwent some evolution and
their properties are now fairly well understood
\cite{Kim_Moin_1985,Rann91,Shen92c,Shen96b,TMV96,MR1440332,Guermond_1999}
(the reader is referred to \cite{MR2250931} for an overview), the same
fundamental idea of decomposing vector fields into a divergence-free
part and a gradient has remained unchanged over the years and has been
challenged only recently in \cite{Guermond20092834}. For all these
schemes, the total cost per time step is that of solving one
vector-valued advection-diffusion equation and one scalar-valued
Poisson equation with homogeneous Neumann boundary conditions.  For
very large size problems, the cost of solving the Poisson equation is
dominant. To address this issue, Guermond and Minev have proposed a
new method in \cite{Guermond2010581}.  The main idea consists of
abandoning the projection paradigm, as in \cite{Guermond20092834}, and
replacing the Poisson equation by a direction splitting strategy. This
requires to solve a sequence of one-dimensional elliptic problems
instead of one multidimensional Poisson equation. The first-order
accurate variant of method has been shown to be unconditionally stable
in \cite{Guermond2010581}.

In this paper we pursue further the ideas introduced/announced in
\cite{Guermond2010581} in the sense that in addition
to splitting the pressure-correction, we also apply a direction
splitting technique to the momentum equation, thus further reducing
the overall computational cost of the method.  We prove that the
totally split method is convergent and we provide error estimates.

Applying direction splitting to the momentum equation is not a new
idea. For instance, in \cite[Section~3.7.2]{MR1846644} Temam studies a
projection method where the solution of the momentum equation is
obtained using direction splitting and the incompressibility
constraint is enforced by means of a Poisson equation. Stability and
convergence of the scheme are proved therein but no error estimates
are provided. Lu, Neittaanm\"aki and Tai show
in~\cite{MR1095644,MR1183413} that this scheme is
$\calO(\tau^{\frac12})$ accurate, $\tau$ being the time-step.  Our
work differs from these previous results mainly in two
directions. First, we adopt a direction splitting strategy for the
computation of the pressure-correction which renders the method
extremely fast and massively parallelizable.  Second, we provide error
estimates for the proposed scheme, and we show that the so-called
standard version of the scheme is $\calO(\tau)$-accurate in all
quantities irrespective of the space dimension and the rotational
version is $\calO(\tau^{\frac32})$-accurate in two space
dimensions. Numerical experiments show that the result holds true also
in three space dimensions and the actual convergence rate is higher
than $\calO(\tau^{\frac32})$ in two and three space dimensions. The
algorithm has been implemented in a parallel code which has been
observed to have optimal weak scalability.  This code has been used to
compute the transient regime on the three-dimensional lid-driven
cavity at $R_e=1000$ and $R_e=5000$ on a mesh composed of $2\,10^9$
nodes on 512 processors only.

This paper is organized as follows. Section~\ref{sub:Notation}
introduces the notation and establishes some preliminary results.  The
new algorithm is described in Section~\ref{sec:Scheme};
two-dimensional and three-dimensional variants of the algorithm are
presented in \S\ref{sub:2d} and \S\ref{sub:3d}, respectively. The
convergence analysis of the standard form of the algorithm is done in
Section~\ref{sec:Error} and the analysis of the rotational form is
done in Section~\ref{sec:ROT}. In Section~\ref{sec:othertimes} we
briefly discuss the BDF2 technique to march in time. Finally, we
present numerical experiments in Section~\ref{sec:NumExp} to
illustrate the performance of this new class of algorithms.

\subsection{Notation and Preliminaries}
\label{sub:Notation}
We consider the time-dependent Stokes system \eqref{eq:Stokes} on the
finite time interval $[0,T]$ and in the cubic domain $\Omega := (0,1)^d$
with $d=2$ or $3$.

We henceforth consider only the time discretization of the system to
simplify the discussion. Handling the space discretization is a
secondary issue, and the reader is referred to
\cite{Guermond_1999,GuermondSalgadoErrorAnalysis} for the techniques
that can be used for this purpose. Let $\tau>0$ be a time step (for
simplicity taken uniform) and let $t_k = k\tau$ for $0\leq k \leq K:=
\lceil T/\tau \rceil$.  Let $E$ be a normed space, with norm
$\|\cdot\|_E.$ For any time-dependent function $\psi :
[0,T]\rightarrow E$, we denote $\psi^k := \psi(t_k)$ and the sequence
$\{\psi^k\}_{k=0,\ldots,K}$ is denoted by $\psi_\tau$. To simplify the
notation we define the time-increment operator $\delta$ by setting
\begin{equation}
  \delta\psi^k := \psi^k - \psi^{k-1},
\end{equation}
and the time-average by
\begin{equation}
  \bar \psi^{k+\frac12} := \frac{ \psi^{k+1} + \psi^k }2.
\end{equation}
We also define the following discrete norms:
\begin{equation}
  \| \psi_\tau\|_{\ell^2(E)} := \left( \tau \sum_{k=0}^K \|\psi^k\|_E^2 \right)^{\frac12},
  \qquad
  \|\psi_\tau\|_{\ell^\infty(E)} := \max_{0\leq k\leq K}\left\{ \|\psi^k\|_E \right\}.
\end{equation}
The space of functions $\psi: [0,T]\longrightarrow E$ that are such
that the map $(0,T)\ni t\longrightarrow \|\psi(t)\|_E\in \Real$ is
$L^p$-integrable is indifferently denoted $L^p((0,T);E)$ or $L^p(E)$.

No notational distinction is done between scalar or vector-valued
functions but spaces of vector-valued functions are identified with
bold fonts. We use the standard Sobolev spaces $W^{m,p}(\Omega),$ for
$0\leq m\leq\infty$ and $1\leq p \leq \infty.$ The closure with
respect to the norm $\|\cdot\|_{W^{m,p}}$ of the space of
$\calC^\infty$-functions compactly supported in $\Omega$ is denoted
$W^{m,p}_0(\Omega)$. To simplify the notation, the Hilbert space
$W^{s,2}(\Omega)$ (resp. $W^{s,2}_0(\Omega)$) is denoted $H^s(\Omega)$
(resp. $H^s_0(\Omega)$). We define $\tildeLdeux$ (resp. $\tildeHun$)
the space that is composed of those functions in $\Ldeux$
(resp. $\Hun$) that are of zero mean.  The scalar product of $\Ldeuxd$
and $\tildeLdeux$ is denoted $\langle\cdot,\cdot\rangle$ and we define
\begin{align}
\|q\|_{H^1} := \|\GRAD q \|_{\bL^2}, \qquad \forall  q\in \tildeHun,\\
\|v\|_{\bH^1}:=\|\GRAD v \|_{\bL^2},\qquad \forall v\in \bH^1_0(\Omega).
\end{align}
Finally we recall that
\begin{equation}
  \| v \|_{\bH^1}^2 = \| \DIV v \|_{L^2}^2 + \| \ROT v \|_{\bL^2}^2,
  \qquad \forall v\in \bH^1_0(\Omega).
\label{eq:grad_eq_divProt}
\end{equation}

Henceforth $c$ denotes a generic constant whose value may change at
each occurrence. This constant may depend on the data of the problem
and its exact solution, but it does not depend on the discretization
parameter $\tau$ or the solution of the numerical scheme.

\subsection{Direction Splitting Pressure Operator} We assume that we have
at hand an operator $A: D(A) \subset \tildeLdeux \rightarrow
\tildeLdeux$ which is unbounded, closed and satisfies
\begin{equation}
  \| \GRAD q \|_{\bL^2}^2 \leq \langle Aq, q \rangle, \quad \forall q \in D(A).
\label{eq:A-coercive}
\end{equation}
This property implies that the map $D(A) \ni q \mapsto \|q\|_A \in \Real$ where
\begin{equation}
  \| q \|_A := \langle Aq,q\rangle^\frac12,\qquad \forall q\in D(A),
\end{equation}
is a norm. We also define the scalar product
\begin{equation}
  \langle p, q\rangle_A := \langle Ap,q\rangle,\qquad \forall p, q \in D(A).
\end{equation}

A natural example for $A$ consists of using $A=-\LAP_N$, where
$\LAP_N$ is the Laplace operator supplemented with homogeneous Neumann
boundary conditions. This operator is the workhorse of classical
projection methods. The main originality of the method that we are going to consider
consists of introducing a direction factorization of this operator.
In two space dimension we define
\begin{equation}
  \begin{cases}
    A := ( 1 - \dxx )( 1 - \dyy ), &  \\
    D(A):= \left\{ p \in \tildeLdeux : \dyy p, Ap \in \Ldeux :
      \partial_y p|_{y=0,1} = 0,\ \partial_x( 1 - \dyy )p|_{x=0,1} = 0
    \right\},
  \end{cases}
\label{eq:A-2d}
\end{equation}
and in three dimensions
\begin{equation}
  \begin{cases}
    A := ( 1 - \dxx )( 1 - \dyy )( 1 - \dzz ), \\
    D(A):= \left\{ p \in \tildeLdeux : \dzz p, ( 1 - \dyy )( 1 - \dzz )p, Ap \in \Ldeux : \right. \\
    \quad \left. \partial_z p|_{z=0,1} =
      0,\ \partial_y(1-\dzz)p|_{y=0,1}
      =0,\ \partial_x(1-\dyy)(1-\dzz)p|_{x=0,1} = 0 \right\},
  \end{cases}
\label{eq:A-3d}
\end{equation}
The graph norm is denoted $\|\cdot\|_{D(A)}$ both in two and three
space dimensions.
\begin{prop}
  The operator $A$ defined in \eqref{eq:A-2d} or \eqref{eq:A-3d}, in
  two or three space dimensions, respectively, satisfies
  \eqref{eq:A-coercive}.
\end{prop}
\begin{proof}
See \cite{Guermond2010581}.
\end{proof}

One interesting feature of the operators defined by \eqref{eq:A-2d}
and \eqref{eq:A-3d} is that solving the equation $Ap=f$ for $f\in
\Ldeux$ only requires to solve one-dimensional problems. For instance,
the solution of $Ap=f$ in three space dimensions is obtained by
solving for $p_1$, $p_2$, and $p$ so that
\begin{align*}
  p_1 - \dxx p_1 &= f,    &\ \partial_x p_1|_{x=0,1} = 0, \\
  p_2 - \dyy p_2 &= p_1,  &\ \partial_y p_2|_{y=0,1} = 0, \\
  p - \dzz p &= p_2, &\ \partial_z p|_{z=0,1} = 0.
\end{align*}

Finally we introduce the Hilbert space $Y$ to be the completion of the
space of smooth scalar-valued functions with respect to the norm
$\|\cdot\|_A$:
\begin{equation}
  Y := \overline{\calC^\infty(\Omega)}^{\|\cdot\|_A}\cap \tildeLdeux.
\end{equation}
The extension of the scalar product $\langle \cdot,\cdot\rangle_A$ to
$Y$ is abusively denoted $\langle \cdot,\cdot\rangle_A$.  For instance
if $A$ is defined as in \eqref{eq:A-2d} or \eqref{eq:A-3d}, the space
$Y$ is characterized as follows:
\begin{equation}
  Y = \begin{cases}
    \{q\in \tildeHun \st \partial_{xy}q\in \Ldeux\},&\text{in }\Real^2,\\
    \{q\in \tildeHun \st \partial_{xy}q, \partial_{yz}q,
    \partial_{zx}q, \partial_{xyz}q \in \Ldeux\}&\text{in }\Real^3.
\end{cases}
\end{equation}
Note that the boundary conditions associated with $D(A)$ have
disappeared from $Y$ and the $\|\cdot\|_A$-norm (which is also the
norm in $Y$) is characterized by
\begin{equation}
  \|q\|_A^2 = \begin{cases}
    \|q\|_{H^1}^2 +\|\partial_{xy}q\|_{L^2}^2, &\text{in }\Real^2,\\
    \|q\|_{H^1}^2 + \|\partial_{xy} q\|_{L^2}^2 + \|\partial_{yz} q\|_{L^2}^2
    + \|\partial_{zx} q\|_{L^2}^2 + \|\partial_{xyz} q\|_{L^2}^2,  &\text{in }\Real^3.
\end{cases}
\end{equation}

\subsection{Direction Splitting Velocity Operator}
To be able to handle the two-dimensional and three-dimensional error
analysis in a unified framework we introduce the following unbounded
closed operator
\begin{equation}
B v :=\begin{cases}
  \partial_{xxyy}v&\text{in $\Real^2$},\\
  (\partial_{xxyy} + \partial_{yyzz} + \partial_{zzxx} -
  \frac{\tau}2 \partial_{xxyyzz})v&\text{in $\Real^3$},
\end{cases} \label{def_of_B}
\end{equation}
with domain
\begin{equation}
D(B) := \{ v\in \Hunzd\st Bv\in \Ldeuxd\}.
\end{equation}
The graph norm is denoted $\|\cdot\|_{D(B)}$.

\begin{lem} The bilinear form $D(B){\times}D(B)\ni (v,w)\longmapsto
  \langle v,Bw\rangle \in \Real$ is symmetric positive and the
  following holds for all $v\in D(B)$:
\begin{equation}
  \langle v,Bv \rangle = \begin{cases}
    \|\partial_{xy}v\|_{\bL^2}^2 &\text{in }\Real^2,\\
    \|\partial_{xy}v\|_{\bL^2}^2 + \|\partial_{yz}v\|_{\bL^2}^2 + \|\partial_{xz}v\|_{\bL^2}^2
    +\frac{\tau}{2} \|\partial_{xyz}v\|_{\bL^2}^2 &\text{in }\Real^3.
\end{cases} \label{def_norm_B}
\end{equation}
\end{lem}
\begin{proof}
  Let us consider the two-dimensional case first. Using the
  Fubini-Tonelli Theorem and integrating by parts repeatedly we obtain
\begin{align*}
  \langle \partial_{xxyy} v, v \rangle &= \int_{x=0}^{x=1} \left[
    v \partial_{yxx} v \Big|_{y=0}^{y=1}
    - \int_{y=0}^{y=1}\partial_y v \partial_{yxx}v  \diff y \right] \diff x \\
  &= -\int_{y=0}^{y=1} \left[ \partial_y v \dyx v \Big|_{x=0}^{x=1}
    - \int_{x=0}^{x=1} (\dyx v)^2  \diff x \right] \diff y \\
  &= \|\dxy v\|_{\bL^2}^2.
\end{align*}
Note that we used $v|_{y=0,1}=0$ and $\partial_y v|_{x=0,1}=0$ which
is a consequence of $v|_{x=0,1}=0$. The three-dimensional result is
obtained similarly; the details are left to the reader.
\end{proof}

To simplify notation we now define the norm
\begin{equation}
  \|v\|_B := \langle v,Bv \rangle^{\frac12}, \qquad v\in D(B),
\end{equation}
and we define the following Hilbert space
\begin{equation}
\bZ := \overline{\calC^\infty(\Omega)}^{\|\cdot\|_B}\cap \Hunzd.
\end{equation}
The extension of the scalar product $\langle \cdot,\cdot\rangle_B$ to
$\bZ$ is abusively denoted $\langle \cdot,\cdot\rangle_B$.

\subsection{The Right-Inverse of the Stokes Operator}
To describe solenoidal vector fields we introduce the classical spaces
\begin{equation}
  \bH := \left\{ v \in \Ldeuxd: \DIV v=0,\ v\SCAL n|_{\partial\Omega} = 0 \right\},
  \quad
  \bV := \bH \cap \Hunzd,
\end{equation}
where $ n$ is the outer unit normal to $\partial\Omega$ and we denote
by $P_\bH$ the $\bL^2$-projection onto $\bH$. It is also useful to
introduce the right-inverse of the Stokes operator $S:\Ldeuxd
\rightarrow \bV$ defined as follows: for any $f\in\Ldeuxd$ we denote $(Sf,q)\in
\bV \times \tildeLdeux$ the pair such that
\begin{equation}
  \begin{cases}
    -\LAP Sf + \GRAD q = f, & \text{in } \Omega, \\
    \DIV Sf =0, & \text{in } \Omega, \\
    Sf = 0, & \text{on } \partial\Omega.
  \end{cases} \label{def_of_right_inverse}
\end{equation}
Given the particular domain that we consider in this work, the inverse
Stokes operator is bounded from $\bL^2(\Omega)$ to $\bH^2\cap \bV$
(\cf \cite{MR977489}), \ie $\|Sf\|_{\bH^2}\le c \|f\|_{\bL^2}$.
Finally, we introduce the seminorm
\begin{equation}
  |v|_\star^2 := \langle Sv,v \rangle, \quad \forall v \in \Ldeuxd,
\end{equation}
and, we recall (see \eg \cite{Shen96b,Guermond_1999,MR2059733}), that for every
$\gamma \in (0,1)$, there exists $c(\gamma)\ge 0$ so that the following
holds for every $v\in \Hunzd$
\begin{equation}
  \langle \GRAD Sv, \GRAD v \rangle \geq  (1-\gamma) \| v \|_{\bL^2}^2
  - c(\gamma)\| v - P_\bH v\|_{\bL^2}^2.
\label{eq:Stokes-prop}
\end{equation}

\section{Description of the Scheme}
\label{sec:Scheme}
We describe the direction splitting algorithm in two and three space
dimensions in this section. The stability and convergence analysis is
done in the subsequent sections.
\subsection{Two Space Dimensions}
\label{sub:2d}
To simplify the presentation, we assume for the time being that the
space dimension is two ($d=2$) and we defer to \S\ref{sub:3d}
the discussion of the three dimensional case.

The scheme computes three sequences of variables $\{u^k\}$,
$\{\phi^{k-\frac12}\}$, and $\{p^{k-\frac12}\}$ that approximate the velocity,
the pressure-correction, and the pressure, respectively.
\begin{itemize}
\item \underline{Pressure predictor:} Denoting by $\pe_0$ the pressure
  field at $t=0$ and $\phi^{\star,-\frac12}$ an approximation of
    $\frac12 \tau \partial_t \pe(0)$, the algorithm is initialized by
    setting $p^{-\frac12}=\pe_0$ and
    $\phi^{-\frac12}=\phi^{\star,-\frac12}$. Then, for all $k\ge 0$ a
    pressure predictor is computed as follows:
\begin{equation}
  p^{\star,k+\frac12} = p^{k-\frac12}+\phi^{k-\frac12}. \label{eq:predictor}
\end{equation}
\item \underline{Velocity update:} The velocity field is initialized
  by setting $u^0=\ue_0$, and for all $k\ge 0$ the velocity update is
  computed by solving the following series of one-dimensional
  problems: Find $u^{k+\frac12}$ and $u^{k+1}$ such that
  \begin{align}
    \frac{u^{k+\frac12} - u^k}{\tau/2} - \dxx u^{k+\frac12} - \dyy u^k + \GRAD
    p^{\star,k+\frac12} = f^{k+\frac12},
    & \quad u^{k+\frac12}|_{x=0,1} = 0, \label{eq:halfstep} \\
    \frac{u^{k+1} - u^{k+\frac12}}{\tau/2} - \dxx u^{k+\frac12} - \dyy u^{k+1}
    + \GRAD p^{\star,k+\frac12} = f^{k+\frac12}, & \quad u^{k+1}|_{y=0,1} =
    0. \label{eq:intstep}
  \end{align}
\item \underline{Penalty step:} The pressure-correction $\phi^{k+\frac12}$
  is computed by solving
  \begin{equation}
    A \phi^{k+\frac12} = -\frac1\tau \DIV u^{k+1}.
    \label{eq:Penalty}
  \end{equation}
\item \underline{Pressure update:} The last sub-step of the algorithm
  consists of updating the pressure as follows:
  \begin{equation}
    p^{k+\frac12} = p^{k-\frac12} + \phi^{k+\frac12} - \chi \DIV \bar u^{k+\frac12}.
  \label{eq:presupdate}
  \end{equation}
\end{itemize}

\begin{rem}
  The parameter $\chi\geq0$ in \eqref{eq:presupdate} is user
  dependent.  By analogy with the projection-based pressure correction
  schemes, we say that the method is in standard form if $\chi=0$
  and the method is in rotational from if $\chi > 0$.
\end{rem}

\begin{rem}
  The splitting of the momentum equation in
  \eqref{eq:halfstep}-\eqref{eq:intstep} is obtained by using the
  original alternating directions (ADI) scheme of Peaceman and
  Rachford, see \cite{MR0071874}.
\end{rem}

\begin{rem}
  The quantity $\phi^{\star,-\frac12}$ can be estimated in many ways.
  For instance one can take $\phi^{\star,-\frac12}=0$; this limits the
  convergence of the scheme to first-order. One can also take
  $\phi^{\star,-\frac12}= \pe^{\star,\frac12} -\pe_0$ where is
  $\pe^{\star,\frac12}$ is an estimate of $\pe(\frac\tau2)$.
\end{rem}

A remarkable feature of the algorithm \eqref{eq:predictor} to
\eqref{eq:presupdate} is that, although the Dirichlet boundary
condition on the velocity is not enforced on the entire boundary at
the integer time steps, it is indeed fulfilled as claimed in the
following
\begin{prop}
\label{prop:bcs}
Let $\{u^k\}$ be the velocity sequence from the algorithm
\eqref{eq:predictor} to \eqref{eq:presupdate}. Then
$u^k|_{\partial\Omega}=0$ for all $k=0,\ldots,K$.
\end{prop}

\begin{proof}
  It is clear that the boundary condition is satisfied at $y=0,1$. Now
  taking the difference of \eqref{eq:intstep} and \eqref{eq:halfstep},
  we obtain the following expression for the half-step velocity
\begin{equation}
  u^{k+\frac12} = \frac{ u^{k+1} + u^k }2 - \frac\tau4\dyy( u^{k+1} - u^k ).
  \label{eq:uhalfeqs}
\end{equation}

Let us consider $x=1$, the other boundary can be treated
similarly. The boundary condition at $x=1$ on the half-step velocity
$u^{k+\frac12}$ implies that
\[
  u^{k+1}(1,y) + u^k(1,y) = \frac\tau2 \dyy ( u^{k+1}(1,y) - u^k(1,y) ).
\]
Moreover, the boundary conditions on $u^{k+1}$ and $u^k$ at $y=1$ imply that
this can be re-written into the following evolution equation
\[
u^{k+1}(1,y) - u^k(1,y) - \frac\tau2 \dyy(u^{k+1}-u^k)(1,y) = -2
u^k(1,y), \quad (u^{k+1}-u^k)(1,\cdot)|_{y=0,1} = 0.
\]
Since $u^0(1,y)=0$ and the evolution operator is positive definite, we
obtain that $u^k(1,y)=0$ for all $k=0,\ldots,K$.
\end{proof}
This result turns out to be crucial for the error analysis.

\subsection{Three Space Dimensions}
\label{sub:3d}
The purpose of this section is to propose a
three-dimensional version of the above splitting technique.  Since the
alternating directions method of Peaceman and Rachford described in
\cite{MR0071874} does not extend to three dimensions, we use the
alternating directions method proposed by Douglas~\cite{MR0136083}
instead to approximate the momentum equation.

The algorithm is again composed of
four steps: pressure predictor, velocity update, penalty step, pressure update.
\begin{itemize}
\item \underline{Pressure predictor:} Denoting by $\pe_0$ the pressure
  field at $t=0$ and $\phi^{\star,-\frac12}$ an approximation of
  $\frac12 \tau \partial_t \pe(0)$, the algorithm is initialized by
  setting $p^{-\frac12}=\pe_0$ and $\phi^{-\frac12}=\phi^{\star,-\frac12}$. Then for all
  $k\ge 0$ a pressure predictor is computed as follows:
\begin{equation}
  p^{\star,k+\frac12} = p^{k-\frac12}+\phi^{k-\frac12}. \label{eq:3dpredictor}
\end{equation}
\item \underline{Velocity update:} The velocity field is initialized
  by setting $u^0=\ue_0$, and for all $k\ge 0$ the velocity update is
  computed by solving the following series of one-dimensional
  problems: Find $\xi^{k+1}$, $\eta^{k+1}$, $\zeta^{k+1}$, and $u^{k+1}$ such that
\begin{align}
  \frac{ \xi^{k+1} - u^k }\tau - \LAP u^k + \GRAD p^{\star,k+\frac12} = f^{k+\frac12}, & \quad \xi^{k+1}|_{\partial\Omega} = 0,
  \label{eq:3dxi} \\
  \frac{ \eta^{k+1} - \xi^{k+1} }\tau - \frac12 \dxx ( \eta^{k+1} - u^k ) = 0, & \quad \eta^{k+1}|_{x=0,1} = 0,
  \label{eq:3deta} \\
  \frac{ \zeta^{k+1} - \eta^{k+1} }\tau - \frac12 \dyy ( \zeta^{k+1} - u^k ) = 0, & \quad \zeta^{k+1}|_{y=0,1} = 0,
  \label{eq:3dzeta} \\
  \frac{ u^{k+1} - \zeta^{k+1} }\tau - \frac12 \dzz ( u^{k+1} - u^k ) = 0, & \quad u^{k+1}|_{z=0,1} = 0.
  \label{eq:3du}
\end{align}
\item \underline{Penalty step:} The pressure-correction $\phi^{k+\frac12}$
  is computed by solving
  \begin{equation}
    A \phi^{k+\frac12} = -\frac1\tau \DIV u^{k+1}.
    \label{eq:3dPenalty}
  \end{equation}
\item \underline{Pressure update:} The last sub-step of the algorithm
  consists of updating the pressure as follows:
  \begin{equation}
    p^{k+\frac12} = p^{k-\frac12} + \phi^{k+\frac12} - \chi \DIV \bar u^{k+\frac12}.
  \label{eq:3dpresupdate}
  \end{equation}
\end{itemize}

\begin{rem} \label{rem:douglas}
  This method is an extension of the alternating direction method
  proposed by Douglas~\cite{MR0136083}. In order to see this, we add
  \eqref{eq:3dxi} and \eqref{eq:3deta} to obtain
\[
  \frac{ \eta^{k+1} - u^k }\tau - \frac12\dxx( \eta^{k+1} + u^k ) - (\dyy+\dzz)u^k + \GRAD p^{\star,k+\frac12} = f^{k+\frac12}.
\]
Adding \eqref{eq:3dxi}--\eqref{eq:3dzeta} we obtain
\[
  \frac{ \zeta^{k+1} - u^k }\tau - \frac12\dxx( \eta^{k+1} + u^k ) - \frac12\dyy( \zeta^{k+1} + u^k )
  -\dzz u^k + \GRAD p^{\star,k+\frac12} = f^{k+\frac12}.
\]
Finally, adding \eqref{eq:3dxi}--\eqref{eq:3du} we obtain
\[
  \frac{ u^{k+1} - u^k }\tau - \frac12\dxx( \eta^{k+1} + u^k ) - \frac12\dyy( \zeta^{k+1} + u^k )
  -\frac12\dzz( u^{k+1} + u^k ) + \GRAD p^{\star,k+\frac12} = f^{k+\frac12}.
\]
These equations correspond to (3.1a)--(3.1c) of \cite{MR0136083}, respectively.
\end{rem}

\begin{prop}\label{prop:3dbcs}
Let $\{u^k\}$ be the velocity sequence from the algorithm
\eqref{eq:3dpredictor}--\eqref{eq:3dpresupdate}. Then
$u^k|_{\partial\Omega}=0$ for all $k=0,\ldots,K$.
\begin{equation}
  u^k|_{\partial\Omega} = 0, \forall k\geq0.
\end{equation}
\end{prop}

\begin{proof}
  Proceed as in the proof of Proposition~\ref{prop:bcs}.
\end{proof}

\subsection{Compatibility Conditions}
Note that $\pe_0:=\pe|_{t=0}$ is not part of the initial data but this
quantity can be computed by solving
\begin{equation}
\LAP \pe_0 = \DIV(f_0 +\LAP \ue_0),\quad
\partial_n \pe_0|_{\front} = (f_0+\LAP \ue_0)\SCAL n,
\end{equation}
where we have set $f_0:=f|_{t=0}$.  This then requires the initial
data to satisfy the following compatibility condition at the boundary
$(- \LAP \ue_0 +\GRAD \pe_0 -f_0)|_{\front} = 0$
which we assume to hold.  This condition holds for instance if $\ue_0=0$
and $f_0=0$, \ie the fluid is a rest at $t=0$ and the source term is
zero at $t=0$.  If the above compatibility condition is not satisfied,
the error analysis must be adapted to account for weighted error
estimates by proceeding as in \cite{MR650052,MR1472237}.

\section{Error Analysis of the Standard Scheme}
\label{sec:Error}
The purpose of this section is to study the convergence of the
algorithms \eqref{eq:halfstep}--\eqref{eq:presupdate} in two space
dimensions and \eqref{eq:3dpredictor}--\eqref{eq:3dpresupdate} in
three space dimensions for $\chi=0$. The main claim of this section
is that the standard version of our scheme is unconditionally stable
and first-order convergent in all quantities.
\subsection{Consistency of the Momentum Equation}
\label{sub:cons_mom}
To evaluate the consistency error on the momentum equation, we
re-write the momentum equation in a more recognizable Crank-Nicolson form.  This is
done in two space dimensions by adding \eqref{eq:halfstep} and
\eqref{eq:intstep} as follows:
\begin{equation}
\frac{u^{k+1} - u^k}{\tau} - \dxx u^{k+\frac12} -\dyy \bar u^{k+\frac12} +
\GRAD p^{\star,k+\frac12} =  f^{k+\frac12}.
\end{equation}
Then using \eqref{eq:uhalfeqs} we obtain the evolution equation for the integer steps,
\begin{equation}
\frac{u^{k+1} - u^k}{\tau} - \LAP \bar u^{k+\frac12} +  \GRAD p^{\star,k+\frac12}
  + \frac{\tau}4\partial_{xxyy} \delta u^{k+1} = f^{k+\frac12}.
\label{eq:discmomentum}
\end{equation}

The same trick can be used in three space dimensions as suggested in
\cite{MR0136083}. By proceeding as above, the intermediate steps,
$\xi^{k+1}$, $\eta^{k+1}$, and $\zeta^{k+1}$, can be eliminated, so
that the momentum equation becomes:
\begin{equation}
  \frac{ u^{k+1} - u^k }\tau - \LAP \bar u^{k+\frac12}
+ \frac\tau4 ( \partial_{xxyy} + \partial_{yyzz} + \partial_{zzxx}
  - \frac{\tau}2 \partial_{xxyyzz})\delta u^{k+1} + \GRAD p^{\star,k+\frac12} = f^{k+\frac12}.
\label{eq:3dmom}
\end{equation}
Owing to the definition of the operator $B$ (see \eqref{def_of_B}),
the momentum equation can be re-written as follows independently of the space dimension:
\begin{equation}
  \frac{ u^{k+1} - u^k }\tau - \LAP \bar u^{k+\frac12}
+ \frac\tau4 B\delta u^{k+1} + \GRAD p^{\star,k+\frac12} = f^{k+\frac12}.
\label{eq:ndmom}
\end{equation}

\subsection{Consistency Analysis of the Algorithm}
\label{sub:cons_rot}
Let $\ue$, $\pe$ be the solution of \eqref{eq:Stokes}.  We define the
following velocity and pressure errors:
\begin{equation}
  e^{k+1} := \ue^{k+1} - u^{k+1}, \quad \epsilon^{k+\frac12} := \pe^{k+\frac12} - p^{k+\frac12},
\end{equation}
where $\ue^{k+1}:=\ue(t_{k+1})$ and
$\pe^{k+\frac12}:=\pe(t_{k+\frac12})$.

Next, we obtain equations controlling the errors.  Since $\chi=0$,
the pressure update implies that the pressure predictor can be written
as follows:
\begin{equation}
\epsilon^{\star,k+\frac12} = 2\epsilon^{k-\frac12}-\epsilon^{k-\frac32}, \qquad \pe^{\star,k+\frac12}
= 2\pe^{k-\frac12}-\pe^{k-\frac32},
\end{equation}
that is, the pressure predictor is a second-order extrapolation of the
pressure at time level $k+\frac12$, and upon subtracting
\eqref{eq:ndmom} from the momentum equation \eqref{eq:Stokes}, we
obtain
\begin{equation}
 (1+  \frac{\tau^2}4 B)(e^{k+1} -e^k) - \tau \LAP \bar e^{k+\frac12} + \tau \GRAD \epsilon^{\star,k+\frac12}
 = \tau \calR^{k+\frac12},
\label{eq:momerr}
\end{equation}
where the residual $\calR^{k+\frac12}$ is defined by
\begin{equation}
  \calR^{k+\frac12} = \left[ \frac{ \delta \ue^{k+1} }\tau - (\ue_t)^{k+\frac12} \right]
            - \LAP\left[ \bar\ue^{k+\frac12} - \ue^{k+\frac12} \right]
            - \GRAD\left[ \pe^{k+\frac12} - \pe^{\star,k+\frac12} \right]
            + \frac\tau4 B \left[ \delta \ue^{k+1} \right].
\end{equation}
Finally, using
\eqref{eq:presupdate} (or \eqref{eq:3dpresupdate}) with $\chi = 0$
to eliminate $\phi^{k+\frac12}$ from \eqref{eq:Penalty} (or
\eqref{eq:3dPenalty}) and using the incompressibility constraint, we
obtain
\begin{equation}
 \langle\delta\epsilon^{k+\frac12},q\rangle_A =
\frac1\tau \langle e^{k+1},\GRAD q\rangle + \langle \delta \pe^{k+\frac12},q\rangle_A\,
\qquad \forall q\in Y.
\label{eq:preserr}
\end{equation}
Note that it is not legitimate to write the equality in strong form,
\ie $A \delta\epsilon^{k+\frac12}$ is equal to $-\frac1\tau \DIV
e^{k+1} + A \delta \pe^{k+\frac12}$ since $\delta\pe^{k+\frac12}$ is
not in $D(A)$ (\ie $\delta\pe^{k+\frac12}$ does not satisfy the
artificial boundary conditions associated with $D(A)$).

\begin{lem}
\label{lem:residual}
Let
$\ue \in W^{2,\infty}(\bH^2(\Omega))\cap W^{1,\infty}(D(B))$
and $\pe \in
W^{2,\infty}(\Hun).$ Then
\begin{equation}
2\tau\langle \calR^{k+\frac12}, v \rangle \leq c\tau^5 + \frac\tau4
\|v\|_{\bL^2}^2, \qquad \forall v\in\Ldeuxd.
\end{equation}
\end{lem}
\begin{proof}
  Each of the terms in $\calR^{k+\frac12}$ is $\calO(\tau^2)$, given
  the smoothness of the exact solution.  Note that $\pe^{k+\frac12} -
  \pe^{\star,k+\frac12} = \delta^2 \pe^{k+\frac12}$.
\end{proof}

\subsection{First Order Estimates on the Velocity}
\label{sub:ErrorFO}
Let us assume that the quantity $\phi^{\star,\frac12}$ is estimated so that
the following holds
\begin{equation}
  \|\pe(\tfrac\tau2) - p^{\star,\frac12}\|_{L^2} \le c \tau. \label{init_pressure_std}
\end{equation}
This is the case if $\phi^{\star,\frac12}=0$ and if the pressure is
smooth enough, say $\pe\in \calC^0([0,T],\Ldeux)$. Then, the main
result of this section is the following first-order convergence
statement:
\begin{thm}
\label{thm:fo}
Assume that the solution $(\ue,\pe)$ to \eqref{eq:Stokes} is smooth
enough (say $\ue \in W^{2,\infty}(\bH^2(\Omega))\cap
W^{1,\infty}(D(B))$ and $\pe\in W^{2,\infty}(Y)$). Then, provided
that \eqref{init_pressure_std} holds, the solution
$(u_\tau,p_\tau)$ to the discrete scheme
\eqref{eq:predictor}--\eqref{eq:presupdate} in two space dimensions
and \eqref{eq:3dpredictor}--\eqref{eq:3dpresupdate} in three space
dimensions, with $\chi=0$, satisfies the following error estimate
\begin{equation}
  \| e_\tau \|_{\ell^\infty(\bL^2)} + \| e_\tau \|_{\ell^2(\bH^1)} + \tau \|e_\tau \|_{\ell^\infty(B)}
  + \tau \| \epsilon_\tau \|_{\ell^\infty(A)} + \sqrt\tau \|\delta e_\tau \|_{\ell^2(B)}
  \leq c \tau.
\end{equation}
\end{thm}
\begin{proof}
  Multiply equation \eqref{eq:momerr} by $2e^{k+1}$ and integrate over
  $\Omega$. Since both the exact velocity and the approximate one at
  integer time steps satisfy the full boundary conditions, we obtain
\begin{multline}
  (1-\frac\tau4)\| e^{k+1} \|_{\bL^2}^2 + \| \delta e^{k+1} \|_{\bL^2}^2
  + \frac\tau2 \| e^{k+1} \|_{\bH^1}^2 + 2\tau \|\bar e^{k+\frac12} \|_{\bH^1}^2
  + 2\tau \langle \GRAD \epsilon^{\star,k+\frac12} , e^{k+1} \rangle +\\
\frac{\tau^2}{4}\| e^{k+1} \|_{B}^2 + \frac{\tau^2}{4} \| \delta e^{k+1} \|_{B}^2
  \leq \|e^k \|_{\bL^2}^2 + \frac\tau2 \| e^k \|_{\bH^1}^2 + \frac{\tau^2}4 \|e^k \|_{B}^2 + c\tau^5.
\label{eq:momfinal}
\end{multline}
Where we have used the identity $2a(a\pm b) = a^2 - b^2 + (a\pm b)^2$.

By using $2\tau^2\epsilon^{\star,k+\frac12}$ as test function in
\eqref{eq:preserr} We obtain
\[
2\tau^2 \langle \delta\epsilon^{k+\frac12},\epsilon^{\star,k+\frac12}
\rangle_A = 2\tau \langle e^{k+1}, \GRAD \epsilon^{\star,k+\frac12}
\rangle + 2\tau^2 \langle \delta \pe^{k+\frac12},
\epsilon^{\star,k+\frac12} \rangle_A.
\]
Clearly,
\[
  \langle \delta\epsilon^{k+\frac12},\epsilon^{\star,k+\frac12} \rangle_A
        = \langle \delta \epsilon^{k+\frac12}, \epsilon^{k+\frac12} \rangle_A
        - \langle \delta \epsilon^{k+\frac12}, \delta^2 \epsilon^{k+\frac12}\rangle_A,
\]
so that using again the identity $2a(a-b) = a^2 - b^2 + (a-b)^2$ we obtain
\begin{multline}
  \tau^2 \left[ \| \epsilon^{k+\frac12} \|_A^2 - \|
    \epsilon^{k-\frac12} \|_A^2 + \| \delta \epsilon^{k-\frac12}
    \|_A^2 \right] - \tau^2 \| \delta^2 \epsilon^{k+\frac12} \|_A^2
 \\ = 2\tau \langle e^{k+1}, \GRAD \epsilon^{\star,k+\frac12} \rangle
  + 2\tau^2 \langle \delta \pe^{k+\frac12}, \epsilon^{\star,k+\frac12}
  \rangle_A.
\label{eq:presfinal}
\end{multline}
To obtain a control on $\| \delta^2 \epsilon^{k+\frac12} \|_A^2$, we
apply the time increment operator $\delta$ to \eqref{eq:preserr}
(assuming that $k\ge 2$) and we use the test function
$\tau\delta^2\epsilon^{k+\frac12}$:
\begin{align*}
  \tau \| \delta^2 \epsilon^{k+\frac12}\|_A^2 &= \langle \delta
  e^{k+1}, \GRAD \delta^2 \epsilon^{k+\frac12} \rangle
  + \tau \langle \delta^2 \pe^{k+\frac12}, \delta^2 \epsilon^{k+\frac12} \rangle_A \\
  & \leq \| \delta e^{k+1} \|_{\bL^2} \|\GRAD \delta^2
  \epsilon^{k+\frac12} \|_{\bL^2}
  + \tau \| \delta^2 \pe^{k+\frac12} \|_A \| \delta^2 \epsilon^{k+\frac12} \|_A \\
  & \leq \left( \| \delta e^{k+1} \|_{\bL^2} + \tau \| \delta^2
    \pe^{k+\frac12}\|_A \right) \| \delta^2 \epsilon^{k+\frac12}\|_A
\end{align*}
So that
\begin{equation}
  \tau^2 \| \delta^2 \epsilon^{k+\frac12} \|_A^2 \leq \| \delta e^{k+1} \|_{\bL^2}^2
  + \tau^2 \| \delta^2 \pe^{k+\frac12} \|_A^2
  + 2\tau \| \delta^2 \pe^{k+\frac12} \|_A \| \delta e^{k+1} \|_{\bL^2}.
\label{eq:shtraf}
\end{equation}

Adding \eqref{eq:momfinal}, \eqref{eq:presfinal} and \eqref{eq:shtraf} we obtain
\begin{multline*}
  (1-\frac\tau4)\| e ^{k+1} \|_{\bL^2}^2 + \frac\tau2 \| e^{k+1}
  \|_{\bH^1}^2 + 2\tau \| \bar e^{k+\frac12} \|_{\bH^1}^2 + \tau^2
  \left[ \| \epsilon^{k+\frac12} \|_A^2
    + \| \delta\epsilon^{k-\frac12} \|_A^2 \right] \\
  + \frac{\tau^2}4 \left[ \| e^{k+1} \|_B^2 + \| \delta e^{k+1} \|_B^2
  \right] \leq c\tau^5 + \| e^k \|_{\bL^2}^2 + \frac\tau2 \| e^k
  \|_{\bH^1}^2 + \tau^2 \| \epsilon^{k-\frac12} \|_A^2
  + \frac{\tau^2}4 \| e^k \|_B^2 \\
  + 2 \tau^2
  \langle\delta\pe^{k+\frac12},\epsilon^{\star,k+\frac12}\rangle_A +
  \tau^2 \| \delta^2 \pe^{k+\frac12} \|_A^2 + 2\tau \| \delta^2
  \pe^{k+\frac12} \|_A \| \delta e^{k+1} \|_{\bL^2}
\end{multline*}
Let us examine the last three terms in detail:\\
$\bullet$ $\tau^2 \| \delta^2 \pe^{k+\frac12} \|_A^2.$ Given the
  smoothness of $\pe$ this term is $\calO(\tau^5)$.\\
$\bullet$ $2\tau \| \delta^2 \pe^{k+\frac12} \|_A \| \delta e^{k+1}
  \|_{\bL^2}.$ We estimate it as follows:
  \[
  2\tau \| \delta^2 \pe^{k+\frac12} \|_A \| \delta e^{k+1} \|_{\bL^2}
  \leq c \tau^3 ( \|e^{k+1}\|_{\bL^2} + \|e^k\|_{\bL^2} ) \leq c
  \tau^5 + \frac\tau4 \|e^{k+1}\|_{\bL^2}^2 + \frac\tau2
  \|e^k\|_{\bL^2}^2.
  \]
$\bullet$ $2 \tau^2 \langle A \delta
  \pe^{k+\frac12},\epsilon^{\star,k+\frac12} \rangle.$ Given the
  smoothness of $\pe$
  \begin{align*}
    2\tau^2\langle\delta\pe^{k+\frac12},\epsilon^{\star,k+\frac12}\rangle_A
    &= 2\tau^2 \langle \delta \pe^{k+\frac12},\epsilon^{k-\frac12}
    \rangle_A
    + 2\tau^2 \langle \delta \pe^{k+\frac12},\delta \epsilon^{k-\frac12} \rangle_A \\
    &\leq
    c\tau^3 \| \epsilon^{k-\frac12} \|_A + c\tau^3 \| \delta \epsilon^{k-\frac12} \|_A \\
    &\leq c\tau^3 + \tau^3 \| \epsilon^{k-\frac12} \|_A^2 + \tau^2 \|
    \delta \epsilon^{k-\frac12} \|_A^2.
  \end{align*}
  Note that this term is the only one in the entire error
  analysis that spoils the game. This consistency term does not allow
  us to obtain directly an error estimate of order larger than
  $\calO(\tau)$.

  We have finally proved that the following holds for all $k\ge 2$:
\begin{multline*}
  (1-\frac\tau2)\| e ^{k+1} \|_{\bL^2}^2 + \frac\tau2 \| e^{k+1}
  \|_{\bH^1}^2 + 2\tau \| \bar e^{k+\frac12} \|_{\bH^1}^2 + \tau^2 \|
  \epsilon^{k+\frac12} \|_A^2
  + \frac{\tau^2}4 \left[ \|e^{k+1}\|_{B}^2 \right. \\
  \left. + \|\delta e^{k+1}\|_{B}^2 \right] \leq c\tau^3 +
  (1+\frac\tau2)\| e^k \|_{\bL^2}^2 + \frac\tau2 \| e^k \|_{\bH^1}^2 +
  \tau^2(1+\tau) \| \epsilon^{k-\frac12} \|_A^2
+ \frac{\tau^2}4\|e^k\|_{B}^2.
\end{multline*}
Upon observing that the initialization process ($p^{-\frac12}=\pe_0$) implies
\[
\tau^2 \|\delta^2 \epsilon^{\frac32}\|_A^2 \le (1+\tau) \|e^2\|^2 + \tau^3,
\]
we infer that the above inequality holds also for $k=1$.  As a
consequence of \eqref{init_pressure_std}, we also deduce that
\begin{align*}
  \|e^{1}\|_{\bL^2}^2 + \tau \|e^{1}\|_{\bH^1}^2 + \tau\|\bar
  e^{\frac12}\|_{\bH^1}^2 +\tau^2 \|\epsilon^{\frac12} \|_A^2 + \tau^2
  \|e^{1}\|_{B}^2 \leq c\tau^4.
\end{align*}
By summing the above relation from $k=1$ to
$K$ and by applying the discrete Gr\"onwall lemma allows us to
conclude.
\end{proof}

The ability of $\delta u^{k+1}/\tau$ to approximate $\ue_t$ is made
explicit in the following:
\begin{lem}
\label{lem:inc}
Let the solution $(\ue,\pe)$ to \eqref{eq:Stokes} be smooth enough
$\ue \in W^{3,\infty}(\bH^2(\Omega))\cap W^{2,\infty}(D(B))$
and $\pe \in W^{3,\infty}(Y)$). Then the solution $(u_\tau,p_\tau)$ to the discrete
scheme \eqref{eq:predictor}--\eqref{eq:presupdate} in two space
dimensions and \eqref{eq:3dpredictor}--\eqref{eq:3dpresupdate} three
space dimensions, with $\chi=0$, satisfies the following error
estimate
\begin{equation}
  \| \delta e_\tau \|_{\ell^\infty(\bL^2)} + \| \delta e_\tau \|_{\ell^2(\bH^1)}
  + \tau \| \delta \epsilon_\tau \|_{\ell^\infty(A)}
  + \tau \|\delta e_\tau\|_{\ell^\infty(B)}
  \leq c \tau^2.
\end{equation}
\end{lem}
\begin{proof}
  Apply the arguments in the proof of Theorem~\ref{thm:fo} to the time
  increments.
\end{proof}

\subsection{Error Estimates on the Pressure}
\label{sub:pres}
It is known that for the incremental projection scheme in standard
form it is possible to prove that the error on the pressure in the
$\ell^2(L^2)$-norm is $\calO(\tau)$ (\cf
\cite{Guermond_1999,MR2250931,Shen96b}). The purpose of this paragraph
is to show that, although on a weaker norm, a similar result holds for
the proposed algorithm. Let us define the norm
\begin{equation}
  \|q\|_{\triangle} =  \sup_{0\neq v \in \bZ } \frac{\langle \GRAD q, v \rangle }{ \| v \|_\bZ}.
\end{equation}
\begin{thm}
  \label{thm:presfoinc}
  Assume that the hypotheses of Lemma~\ref{lem:inc} hold, then
\begin{equation}
\| \epsilon_\tau \|_{\ell^2(\triangle)} \leq c \tau.
\end{equation}
\end{thm}
\begin{proof}
  Using the error equation \eqref{eq:momerr} we obtain
\begin{align*}
  \| \epsilon^{\star,k+\frac12} \|_{\triangle} &= \sup_{0\neq v \in \bZ } \frac1{
    \| v \|_\bZ } \left[ \langle \frac{\delta e^{k+1}}\tau, v \rangle
    + \langle \GRAD \bar e^{k+\frac12}, \GRAD v \rangle + \frac\tau4
    \langle\delta e^{k+1},v\rangle_B + \langle
    \calR^{k+\frac12}, v \rangle
  \right] \\
  & \leq \frac{ \| \delta e^{k+1} \|_{\bL^2} }\tau + \| \bar
  e^{k+\frac12} \|_{\bH^1}
  + \frac\tau4 \|\delta e^{k+1} \|_{B} + c\tau^2 \\
  & \leq c \tau + \| \bar e^{k+\frac12} \|_{\bH^1},
\end{align*}
where the last estimate holds in view of Lemma~\ref{lem:inc}. Take the
square of this inequality, multiply it by $\tau$ and sum over $k$.
The result follows by using the conclusion of Theorem~\ref{thm:fo}.
\end{proof}

\begin{rem}
  It seems that it may be possible to obtain a first-order error
  estimate on the pressure in the $\ell^2(L^2)$-norm in the fully
  discrete case under the additional (somewhat restrictive) condition
\begin{equation}
\tau \le
\begin{cases}
c h & \text{in $\Real^2$},\\
c h^{\frac43} & \text{in $\Real^3$}.
\end{cases} \label{cfl_2d_3d}
\end{equation}
This is a CFL condition in two space dimensions.  The reasoning behind
this conjecture is the following. Assume that the velocity is
approximated using a finite-dimensional space $\bX_h$ and that the
norm in $B$ is appropriately approximated, say $\|\cdot\|_{B_h}$.
In view of \eqref{def_norm_B} it is reasonable to expect that the
following inverse inequalities hold:
\[
\|v\|_{B_h} \le
\begin{cases}
c h^{-1} \|v\|_{\bH^1} & \text{in $\Real^2$},\\
 c h^{-1}(1+ \tau^{\frac12} h^{-1}) \|v\|_{\bH^1} & \text{in $\Real^3$},
\end{cases} \qquad \forall v \in \bX_h.
\]
Then, assuming that the pressure is approximated using a space $M_h
\subset \tildeHun$ so that the pair $(\bX_h,M_h)$ satisfies the
so-called LBB condition, \cite{MR2050138,GR86}, we obtain
\begin{align*}
  c \| \epsilon^{\star,k+\frac12} \|_{L^2}
  &\leq \sup_{0 \neq v \in \bX_h } \frac{ \langle \GRAD \epsilon^{\star,k+\frac12}, v \rangle }{ \| v \|_{\bH^1} } \\
  &\leq \sup_{0 \neq v \in \bX_h } \frac1{ \| v \|_{\bH^1} } \left[ \langle
    \delta e^{k+1}/\tau, v \rangle + \langle \GRAD \bar e^{k+\frac12},
    \GRAD v \rangle + \frac\tau4 \langle\delta e^{k+1},v \rangle_{B_h} +
    \langle \calR^{k+\frac12}, v \rangle
  \right] \\
  &\leq \| \delta e^{k+1}/\tau \|_{\bL^2} + \| \bar
  e^{k+\frac12} \|_{\bH^1} + \frac\tau4 \sup_{0 \neq v \in \bX_h} \frac{
    \|\delta e^{k+1}\|_{B_h} \|v\|_{B_h} }{ \| v \|_{\bH^1} } + c\tau^2.
\end{align*}
The two-dimensional inverse inequality implies
\[
\| \epsilon^{\star,k+\frac12} \|_{L^2} \leq c\tau^2 + \|\tau^{-1} \delta
e^{k+1}\|_{\bL^2} +\| \bar e^{k+\frac12} \|_{\bH^1} + c\tau^2 h^{-2} \|\tau^{-1}\delta e^{k+1}\|_{\bH^1},
\]
whereas the  three-dimensional inverse inequality implies
\[
\| \epsilon^{\star,k+\frac12} \|_{L^2} \leq c\tau^2 + \|\tau^{-1} \delta
e^{k+1}\|_{\bL^2} +\| \bar e^{k+\frac12} \|_{\bH^1}
+ c(\tau^2 h^{-2}+\tau^3 h^{-4})\|\tau^{-1}\delta e^{k+1}\|_{\bH^1}.
\]
Take the square of this inequality, multiply it by $\tau$ and sum over
$k$, then the estimates of Lemma~\ref{lem:inc} together with condition
\eqref{cfl_2d_3d} yield the desired estimate,
$\|\epsilon_\tau\|_{\ell^2(L^2)} \le c\tau$.
\end{rem}

\subsection{Second-Order Estimates on the Velocity}
\label{sub:velL2std}
Despite the fact that numerical experiments suggest that the standard
form of the above algorithm is close to second-order on the velocity
in the $\bL^2$-norm, (see Section~\ref{sec:NumExp}), a proof of such
statement eludes us at the moment. We briefly elaborate in this
section on the difficulties that arise when trying to establish a
second-order error estimate.

The argument one usually invokes to prove a second-order error
estimate consists of multiplying the error equation by $S\bar
e^{k+\frac12},$ where $S$ is the right-inverse Stokes operator (see
\eqref{def_of_right_inverse}). Following this reasoning, and using
property \eqref{eq:Stokes-prop}, we obtain that the following holds
\begin{multline}
  \frac12 \left( | e^{k+1} |_\star^2 - |e^k|_\star^2 \right) + \frac{3\tau}{4}\| \bar e^{k+\frac12} \|_{\bL^2}^2
  + \frac{\tau^2}4 \langle \delta e^{k+1}, S\bar e^{k+\frac12} \rangle_B
  \leq \tau \langle \calR^{k+\frac12}, S\bar e^{k+1}\rangle \\
  + c\tau \| \bar e^{k+\frac12} - P_\bH \bar e^{k+\frac12} \|_{\bL^2}^2,
\label{eq:errorxStokes}
\end{multline}
Provided the exact solution is smooth enough, we can estimate the residual term
in a way similar to Lemma~\ref{lem:residual},
\[
\tau\langle\calR^{k+\frac12}, S\bar e^{k+\frac12}\rangle\leq c \tau^5
+ \frac\tau8 \|\bar e^{k+\frac12}\|_{\bL^2}^2.
\]
Using the estimates of Lemma~\ref{lem:inc} we can control the $B$-norm
as follows:
\begin{align*}
 \left |\frac{\tau^2}4 \langle\delta e^{k+1},S \bar e^{k+\frac12} \rangle_B \right|
 &\leq \frac{\tau^2}4 \|\delta e^{k+1}\|_{B} \| S \bar e^{k+\frac12} \|_{B}
 \leq c \tau^3 \| S \bar e^{k+\frac12} \|_{B}.
\end{align*}
In two space dimensions the $\bH^2$-regularity of $S$ implies
$\|S\bar e^{k+\frac12}\|_{B} \le c \|\bar e^{k+\frac12} \|_{\bL^2}$
so that
\begin{align*}
  \left |\frac{\tau^2}4 \langle\delta e^{k+1},S \bar e^{k+\frac12}
    \rangle_B \right| & \leq c \tau^5 + \frac{\tau}{8}\|\bar
  e^{k+\frac12} \|_{\bL^2}^2.
\end{align*}
Note that the above reasoning does not apply in three space
dimensions.  In conclusion, in two space dimensions
\eqref{eq:errorxStokes} becomes
\[
  | e^{k+1} |_\star^2 - |e^k|_\star^2 + \tau \| \bar e^{k+\frac12} \|_{\bL^2}^2
  \leq c\tau(\tau^4 +  \inf_{v \in \bH}\| \bar e^{k+\frac12} - v \|_{\bL^2}^2)
  = c\tau(\tau^4 + \| \bar e^{k+\frac12} - P_\bH \bar e^{k+\frac12}\|_{\bL^2}^2).
\]
which in turn yields
\begin{equation}
\| \bar e_\tau \|_{\ell^2(\bL^2)}^2 \le  c(\tau^4 + \|\bar e_\tau-P_\bH\bar e_\tau\|_{\ell^2(\bL^2)}^2).
\label{eq:L2est}
\end{equation}
This inequality shows that the estimate on $\| \bar e_\tau
\|_{\ell^2(\bL^2)}$ is controlled by $\|\bar e_\tau-P_\bH\bar
e_\tau\|_{\ell^2(\bL^2)}$. Let us now try to bound $\|\bar
e_\tau-P_\bH\bar e_\tau\|_{\ell^2(\bL^2)}$ uniformly.

By definition, there is $\mu^{k+\frac12}\in \tildeHun$ so that $\bar
e^{k+\frac12} - P_\bH \bar e^{k+\frac12} = \GRAD \mu^{k+\frac12}$. In
other words $\mu^{k+\frac12}$ solves $- \LAP \mu^{k+\frac12} = - \DIV
\bar e^{k+\frac12}$ and $\partial_n \mu^{k+\frac12}|_{\partial\Omega}
= 0$. Then the penalty equation \eqref{eq:preserr} together with the
assumed smoothness of the pressure and the estimates of
Lemma~\ref{lem:inc} imply that
\begin{align*}
  \|\bar e^{k+\frac12} - P_\bH \bar e^{k+\frac12}\|_{\bL^2}^2
  &= \|\GRAD\mu^{k+\frac12}\|_{\bL^2}^2 \\
  & = \langle \bar e^{k+\frac12},\GRAD\mu^{k+\frac12}\rangle =
  \frac\tau2 \langle \delta \epsilon^{k+\frac12} - \delta
  \pe^{k+\frac12} +
  \delta\epsilon^{k-\frac12} -  \delta\pe^{k-\frac12},\mu^{k+\frac12}\rangle_A\\
  &\leq \frac\tau2 \left( \| \delta \epsilon^{k+\frac12} \|_A + \|
    \delta \pe^{k+\frac12} \|_A
    + \| \delta \epsilon^{k-\frac12} \|_A + \| \delta \pe^{k-\frac12} \|_A \right) \| \mu^{k+\frac12} \|_A \\
  & \leq c\tau^2 \| \mu^{k+\frac12} \|_A.
\end{align*}
This finally gives the estimate
\[
  \| \bar e^{k+\frac12} - P_\bH \bar e^{k+\frac12} \|_{\bL^2}
      \leq c \tau^2 \frac{ \| \mu^{k+\frac12} \|_A }{ \| \GRAD \mu^{k+\frac12} \|_{\bL^2} },
\]
which can be controlled uniformly if $\|\cdot\|_A$ induces a norm
equivalent to $H^1$. This is unfortunately not true with the operators
$A$ defined in \eqref{eq:A-2d} and \eqref{eq:A-3d}.

In conclusion, the reasoning carried out above seems to indicate that
the right-inverse Stokes operator $S$ is not the correct operator that
should be used for the duality argument. The operator that should be
used instead still eludes us at the moment.

\section{Error Analysis of the Rotational Scheme}
\label{sec:ROT}
The purpose of this section is to analyze the algorithms
\eqref{eq:predictor}--\eqref{eq:presupdate} and
\eqref{eq:3dpredictor}--\eqref{eq:3dpresupdate} for $\chi \neq 0$
and to show that, as it is the case for the classical rotational
pressure-correction schemes (\cf \cite{MR2059733}), these algorithms
provide a better order of convergence than the standard form.

\subsection{Consistency Analysis}
\label{sub:rot_cons}
Let $\ue$, $\pe$ be the solution of \eqref{eq:Stokes}.  We define the
following velocity and pressure errors:
\begin{equation}
  e^{k+1} := \ue^{k+1} - u^{k+1}, \quad
  \epsilon^{k+\frac12} := \pe^{k+\frac12} - p^{k+\frac12},
\end{equation}
where $\ue^{k+1}:=\ue(t_{k+1})$ and
$\pe^{k+\frac12}:=\pe(t_{k+\frac12})$. The error on the pressure correction is
measured by introducing the following quantity:
\begin{equation}
  \calP^{k+\frac12}:=\delta \epsilon^{k+\frac12} + \chi \DIV \bar e^{k+\frac12}.
\end{equation}

Using the above notation we infer
\begin{align*}
  \pe^{k+\frac12} - p^{\star,k+\frac12} &= \pe^{k+\frac12} -
  (p^{k-\frac12} + \phi^{k-\frac12}) = \pe^{k+\frac12} -
  (p^{k-\frac12} + p^{k-\frac12} - p^{k-\frac32}
  + \chi \DIV\bar u^{k-\frac12})\\
  &=\delta^2\pe^{k+\frac12} + (\epsilon^{k-\frac12}+\delta
  \epsilon^{k-\frac12}
  + \chi \DIV\bar e^{k-\frac12})\\
  &= \delta^2\pe^{k+\frac12} + \epsilon^{k-\frac12}+ \calP^{k-\frac12}
\end{align*}
Then momentum equation is rewritten as
follows:
\begin{equation}
  (1+  \frac{\tau^2}4 B)(e^{k+1} -e^k) - \tau \LAP \bar e^{k+\frac12}
  + \tau \GRAD(\epsilon^{k-\frac12}+ \calP^{k-\frac12})
  = \tau \calR^{k+\frac12},
\label{eq:rot_momerr}
\end{equation}
where the residual $\calR^{k+\frac12}$ is defined by
\begin{equation}
  \calR^{k+\frac12} = \left[ \frac{ \delta \ue^{k+1} }\tau - (\ue_t)^{k+\frac12} \right]
            - \LAP\left[ \bar\ue^{k+\frac12} - \ue^{k+\frac12} \right]
            - \GRAD\left[ \delta^2\pe^{k+\frac12} \right]
            + \frac\tau4 B \left[ \delta \ue^{k+1} \right].
\end{equation}
The result of
Lemma~\ref{lem:residual} holds again, $\calR^{k+\frac12} =
\calO(\tau^2)$, provided the exact solution is smooth enough.  The
equation that controls the pressure correction is re-written as
follows:
\begin{equation}
 \langle\calP^{k+\frac12} ,q\rangle_A =
\frac1\tau \langle e^{k+1},\GRAD q\rangle + \langle \delta \pe^{k+\frac12},q\rangle_A\,
\qquad \forall q\in Y.
\label{eq:rot_preserr}
\end{equation}

\subsection{A Priori Estimate on the Divergence of the Velocity}
\label{sub:DIV}
Let us assume the that quantity $\phi^{\star,\frac12}$ is estimated so that
the following holds
\begin{equation}
  \|\pe(\tfrac\tau2) - p^{\star,\frac12}\|_{L^2} \le c \tau. \label{init_pressure_rot}
\end{equation}
This is the case if $\phi^{\star,\frac12}=0$ and if the pressure is
smooth enough, say $\pe\in \calC^0([0,T],\Ldeux)$.
The main result of this section is the
following
\begin{thm}
\label{thm:DIV-est}
Assume that the solution $(\ue,\pe)$ to \eqref{eq:Stokes} is smooth
enough, (say $\ue \in W^{2,\infty}(\bH^2(\Omega))\cap
W^{1,\infty}(D(B))$ and $\pe \in W^{2,\infty}(Y)$). Then, provided
\eqref{init_pressure_rot} holds, the solution $(u_\tau,p_\tau)$ to the
discrete scheme \eqref{eq:predictor}--\eqref{eq:presupdate} in two
space dimensions and \eqref{eq:3dpredictor}--\eqref{eq:3dpresupdate}
in three space dimensions, with $0 < \chi \leq 1$, satisfies the
following error estimate
\begin{align}
  \|\delta e_\tau \|_{{\ell^\infty(L^2)}}^2 + \tau \| \ROT \delta e_\tau
  \|_{\ell^\infty(\bL^2)}^2 + \|\ROT\delta \overline{e}_\tau
  \|_{\ell^2(\bL^2)}^2 + \tau \|\DIV e_\tau \|_{\ell^\infty(L^2)}^2
  \leq c \tau^{4}. \label{eq:thm:DIV-est}
%\\ \|\ue_\tau-u_\tau\|_{{\ell^\infty(L^2)}} +
  %\|\ue_\tau-u_\tau\|_{\ell^\infty(\bH^1)} \leq c \tau^{\frac32}.
\end{align}
\end{thm}
\begin{proof}
  Following \cite{MR2059733}, we derive an improved estimate on the
  divergence of the velocity.  This is done by working with the time
  increments of \eqref{eq:rot_momerr}-\eqref{eq:rot_preserr}.

  Apply the time increment operator $\delta$ to the momentum equation
  \eqref{eq:rot_momerr} and test against $2\delta e^{k+1}$ to obtain
\begin{multline}
  \left(1-\frac\tau4\right)\| \delta e^{k+1} \|_{\bL^2}^2 + \|
  \delta^2 e^{k+1} \|_{\bL^2}^2 + \frac\tau2 \| \delta e^{k+1}
  \|_{\bH^1}^2 + 2\tau\|\delta\bar e^{k+\frac12}\|_{\bH^1}^2
  +\frac{\tau}{4}\|\delta e^{k+1}\|_B^2 \\
  +2\tau \left\langle \GRAD \left( \delta\epsilon^{k-\frac12} + \delta
      \calP^{k-\frac12} \right),\delta e^{k+1} \right\rangle \le
  c\tau^5 + \|\delta e^k \|_{\bL^2}^2 + \frac\tau2 \| \delta e^k
  \|_{\bH^1}^2 +\frac{\tau}{4}\|\delta e^{k}\|_B^2,
\label{eq:rot_dos}
\end{multline}
where we used the fact that the residual is $\calO(\tau^2)$. Note that
we could decrease the consistency error to $\calO(\tau^3)$ by assuming
more regularity on $\ue$ and $\pe$, but it would not improve the
overall accuracy of the method since the splitting error will turn out
to be $\calO(\tau^2)$ (see below).

Apply the time increment operator $\delta$ to \eqref{eq:rot_preserr}
and use $2\tau^2\calP^{k+\frac12}$ as a test function. We obtain
\begin{align*}
  \tau^2 \| \calP^{k+\frac12} \|_A^2 + \tau^2 \| \delta
  \calP^{k+\frac12} \|_A^2 - \tau^2 \| \calP^{k-\frac12} \|_A^2 &=
  -2\tau \langle \DIV \delta e^{k+1}, \calP^{k+\frac12}\rangle
  + 2\tau^2 \langle \delta^2 \pe^{k+\frac12}, \calP^{k+\frac12} \rangle_A \\
  &= -2\tau \langle \DIV \delta e^{k+1}, \chi\DIV\bar
  e^{k+\frac12} + \delta \epsilon^{k+\frac12}\rangle + 2\tau^2 \langle
  \delta^2 \pe^{k+\frac12}, \calP^{k+\frac12} \rangle_A,
\end{align*}
where we used the identity $2a(a \pm b) = a^2 + (a \pm b)^2
-b^2$. This gives
\begin{multline}
  \tau^2 \| \calP^{k+\frac12} \|_A^2 + \tau^2 \| \delta \calP^{k+\frac12} \|_A^2
  + \chi\tau \| \DIV e^{k+1} \|_{L^2}^2
  =
  \tau^2 \| \calP^{k-\frac12} \|_A^2 \\
  + \chi\tau\|\DIV e^k\|_{L^2}^2 +
  2\tau \langle \delta e^{k+1}, \GRAD \delta \epsilon^{k+\frac12} \rangle
  + 2\tau^2 \langle \delta^2 \pe^{k+\frac12}, \calP^{k+\frac12} \rangle_A.
\label{eq:rot-uno}
\end{multline}

Apply again the time increment operator $\delta$ to
\eqref{eq:rot_preserr} and test with
$-2\tau^2\delta^2\calP^{k+\frac12}$. Using, again the identity
$2a(a-b) = a^2 +(a-b)^2 - b^2$, we obtain
\begin{multline}
  -\tau^2 \left[ \| \delta \calP^{k+\frac12}\|_A^2 + \| \delta^2
    \calP^{k+\frac12}\|_A^2 - \| \delta \calP^{k-\frac12}\|_A^2\right]
  = -2\tau \langle \GRAD \delta^2 \calP^{k+\frac12}, \delta e^{k+1}
  \rangle
  \\
  -2\tau^2 \langle \delta^2 \pe^{k+\frac12}, \delta^2
  \calP^{k+\frac12} \rangle_A.
\label{eq:rot_tres}
\end{multline}
Observe that \eqref{eq:rot-uno}+\eqref{eq:rot_tres} amounts to testing
the time increment of \eqref{eq:rot_preserr} with
$2\tau^2(\calP^{k-\frac12} + \delta \calP^{k-\frac12})$. We have split
the two steps to make the argument clearer.

By summing \eqref{eq:rot_dos}, \eqref{eq:rot-uno} and
\eqref{eq:rot_tres} we deduce that
\begin{multline}
  \left(1-\frac\tau4 \right)\| \delta e^{k+1} \|_{\bL^2}^2 +
  \frac\tau2 \| \delta e^{k+1} \|_{\bH^1}^2 + 2\tau \| \delta \bar
  e^{k+\frac12} \|_{\bH^1}^2 +\frac{\tau}{4}\|\delta e^{k+1}\|_B^2
  + \chi \tau \| \DIV e^{k+1} \|_{L^2}^2 \\
  + \tau^2 \| \calP^{k+\frac12} \|_A^2 + \tau^2 \| \delta
  \calP^{k-\frac12} \|_A^2 + \| \delta^2 e^{k+1} \|_{\bL^2}^2 - \tau^2
  \| \delta^2 \calP^{k+\frac12}\|_A^2 - 2\chi\tau \langle \DIV
  \delta \bar e^{k+\frac12}, \DIV \delta e^{k+1} \rangle
  \leq  c\tau^5 \\
  +\| \delta e^k \|_{\bL^2}^2 +\frac\tau2 \| \delta e^k \|_{\bH^1}^2
  +\frac{\tau}{4}\|\delta e^k\|_B^2 + \chi \tau \| \DIV e^k
  \|_{L^2}^2 + \tau^2 \| \calP^{k-\frac12} \|_A^2 + 2 \tau^2 \langle
  \delta^2 \pe^{k+\frac12}, \calP^{k-\frac12} + \delta
  \calP^{k-\frac12} \rangle_A.
\label{eq:rot_cuatro}
\end{multline}
where we used the following identities
\begin{align*}
  \epsilon^{k-\frac12} + \delta \calP^{k-\frac12} -\delta
  \epsilon^{k+\frac12} + \delta^2 \calP^{k+\frac12} = \chi \DIV
  \delta
  \bar e^{k+\frac12}, \\
  \calP^{k+\frac12} - \delta^2 \calP^{k+\frac12} = \calP^{k-\frac12} +
  \delta \calP^{k-\frac12}.
\end{align*}
Given the smoothness of $\pe$, the following holds:
\[
2\tau^2 \langle \delta^2 \pe^{k+\frac12}, \calP^{k-\frac12} + \delta
\calP^{k-\frac12} \rangle_A \leq c\tau^5 + \frac{\tau^3}2 \|
\calP^{k-\frac12} \|_A^2 + \tau^2 \|\delta\calP^{k-\frac12}\|_A^2.
\]
Observe that it is here that the irreducible splitting error comes
into full light.  Although the consistency of the time increment of
the momentum equation is $\calO(\tau^3)$ (provided enough regularity
is assumed on $\ue$ and $\pe$), the above inequality shows that
splitting error of the method is $\calO(\tau^2)$. Then
\eqref{eq:rot_cuatro} becomes
\begin{multline}
  \left(1-\frac\tau4 \right)\| \delta e^{k+1} \|_{\bL^2}^2 +
  \frac\tau2 \| \delta e^{k+1} \|_{\bH^1}^2
- \chi\frac\tau2 \| \DIV \delta e^{k+1} \|_{L^2}^2
+ 2\tau \| \delta \bar
  e^{k+\frac12} \|_{\bH^1}^2
- 2\tau\chi \| \DIV \delta
  \bar e^{k+\frac12} \|_{L^2}^2 \\
+\frac{\tau}{4}\|\delta e^{k+1}\|_B^2
  + \chi \tau \| \DIV e^{k+1} \|_{L^2}^2
  + \tau^2 \| \calP^{k+\frac12} \|_A^2  + \| \delta^2 e^{k+1} \|_{\bL^2}^2
- \tau^2\| \delta^2 \calP^{k+\frac12}\|_A^2
  \leq  c\tau^5 \\
  +\| \delta e^k \|_{\bL^2}^2 +\frac\tau2 \| \delta e^k \|_{\bH^1}^2
-\frac\tau2\chi \| \DIV \delta e^k \|_{L^2}^2
  +\frac{\tau}{4}\|\delta e^k\|_B^2 + \chi \tau \| \DIV e^k
  \|_{L^2}^2 + \tau^2(1+\frac{\tau}{2})\| \calP^{k-\frac12} \|_A^2.
\label{eq:rot_five}
\end{multline}
where we used
\[
-2\chi \langle \DIV \delta \bar e^{k+\frac12}, \DIV \delta e^{k+1}
\rangle = -\chi \left[ \frac\tau2 \| \DIV \delta e^{k+1} \|_{L^2}^2
  - \frac\tau2 \| \DIV \delta e^k \|_{L^2}^2 + 2\tau \| \DIV \delta
  \bar e^{k+\frac12} \|_{L^2}^2 \right].
\]
Then the identity \eqref{eq:grad_eq_divProt} gives
\begin{multline}
  \left(1-\frac\tau4 \right)\| \delta e^{k+1} \|_{\bL^2}^2 +
  \frac\tau2 \|\ROT \delta e^{k+1} \|_{\bL^2}^2 +(1-\chi)\frac\tau2
  \|\DIV\delta e^{k+1} \|_{L^2}^2 + 2\tau \|\ROT\delta \bar
  e^{k+\frac12} \|_{\bH^1}^2 \\ + 2\tau(1-\chi) \| \DIV \delta
  \bar e^{k+\frac12} \|_{L^2}^2
  +\frac{\tau}{4}\|\delta e^{k+1}\|_B^2 + \chi \tau \| \DIV e^{k+1}
  \|_{L^2}^2 + \tau^2 \| \calP^{k+\frac12} \|_A^2 + \| \delta^2
  e^{k+1} \|_{\bL^2}^2 - \tau^2\| \delta^2 \calP^{k+\frac12}\|_A^2
  \leq  c\tau^5 \\
  +\| \delta e^k \|_{\bL^2}^2 + \frac\tau2 \|\ROT\delta e^k
  \|_{\bL^2}^2 +(1-\chi)\frac\tau2 \|\DIV \delta e^k\|_{L^2}^2
  +\frac{\tau}{4}\|\delta e^k\|_B^2 + \chi \tau \| \DIV e^k
  \|_{L^2}^2 + \tau^2(1+\frac{\tau}{2})\| \calP^{k-\frac12} \|_A^2.
\label{eq:rot_six}
\end{multline}

To conclude we are going to observe that the quantity $\| \delta^2
e^{k+1} \|_{\bL^2}^2 - \tau^2 \| \delta^2 \calP^{k+\frac12}\|_A^2$ is
non negative up to some consistency error.  To see this, let us apply
the time increment operator $\delta^2$ to \eqref{eq:rot_preserr} and
test the equation with $\tau\delta^2\calP^{k+\frac12}$. After using
the Cauchy-Schwarz inequality and the inequality
\eqref{eq:A-coercive}, we obtain
\[
\tau \| \delta^2\calP^{k+\frac12} \|_A \leq \| \delta^2 e^{k+1}
\|_{\bL^2} + \tau \| \delta^3 \pe^{k+\frac12} \|_A,
\]
which, given the smoothness assumption on $\pe$, then implies
\begin{align*}
\tau^2 \| \delta^2\calP^{k+\frac12} \|_A^2
&\leq \| \delta^2 e^{k+1}
\|_{\bL^2}^2 + \tau^2 \| \delta^3 \pe^{k+\frac12} \|_A^2 + 2\tau
\|\delta^2 e^{k+1} \|_{\bL^2}\|\delta^3 \pe^{k+\frac12}\|_A,\\
 & \leq c\tau^5 + \| \delta^2 e^{k+1} \|_{\bL^2}^2
  + \frac\tau4 \| \delta e^{k+1} \|_{\bL^2}^2 + \frac\tau2 \| \delta e^k \|_{\bL^2}^2.
\end{align*}
Note again that the consistency error could be decreased to
$\calO(\tau^3)$ by assuming $\pe\in W^{3,\infty}(Y)$, but this would be
useless since the splitting error of the method has been shown to be
$\calO(\tau^2)$ above.
By adding this last inequality to \eqref{eq:rot_six} we finally obtain
that the following holds for all $k\ge 2$.
\begin{multline*}
  \left(1-\frac\tau2 \right)\| \delta e^{k+1} \|_{\bL^2}^2 +
  \frac\tau2 \|\ROT\delta e^{k+1} \|_{\bL^2}^2
 +(1-\chi)\frac\tau2 \|\DIV\delta e^{k+1} \|_{L^2}^2
+ 2\tau \| \ROT \delta \bar
  e^{k+\frac12} \|_{\bL^2}^2 + 2\tau(1-\chi) \| \DIV \delta
  \bar e^{k+\frac12} \|_{L^2}^2 \\
+\frac{\tau}{4}\|\delta e^{k+1}\|_B^2
  + \chi \tau \| \DIV e^{k+1} \|_{L^2}^2
  + \tau^2 \| \calP^{k+\frac12} \|_A^2
  \leq  c\tau^5
  + (1+\frac{\tau}{2})\| \delta e^k \|_{\bL^2}^2  + \frac\tau2 \|\ROT\delta e^k \|_{\bL^2}^2\\
+(1-\chi)\frac\tau2 \|\DIV \delta e^k\|_{L^2}^2
  +\frac{\tau}{4}\|\delta e^k\|_B^2 + \chi \tau \| \DIV e^k
  \|_{L^2}^2 + \tau^2(1+\frac{\tau}{2})\| \calP^{k-\frac12} \|_A^2.
\end{multline*}
The following estimate also holds as a consequence of the
initialization hypothesis \eqref{init_pressure_rot}:
\[
\|\delta e^2\|_{\bL^2}^2  + \tau \|\ROT\delta e^2\|_{\bL^2}^2\\
+\tau \|\DIV \delta e^2\|_{L^2}^2 +\tau\|\delta e^2\|_B^2 + \chi
\tau \| \DIV e^2 \|_{L^2}^2 + \tau^2\|\calP^{\frac32}\|_A^2 \le
\tau^4.
\]

By summing the above inequalities from $k=2$ to $K$ and by applying
the discrete Gr\"onwall lemma we finally obtain the following error
bound:
%\begin{multline*}
\[
\| \delta e_\tau \|_{\ell^\infty(\bL^2)}^2 + \tau \|\ROT
\delta e_\tau\|_{\ell^\infty(\bL^2)}^2 +  \tau \|\ROT \delta \bar
e_\tau\|_{\ell^2(\bL^2)}^2 +
\chi \tau \| \DIV e_\tau\|_{\ell^\infty(L^2)}^2 + \tau^2 \| \calP_\tau
\|_{\ell^\infty(A)}^2 \leq c\tau^4
\]
%\end{multline*}
This completes the proof.
\end{proof}

\subsection{Error Estimates}
\label{sub:rot_Err}
Having obtained the estimate of Theorem~\ref{thm:DIV-est} we can now
show that the rotational version of the algorithm provides a better
order of convergence for the velocity in the $\ell^2(\bL^2)$-norm, at
least in two space dimensions. To this end, let us denote by $\bar
\psi_\tau$ the sequence whose generic term is
$\bar\psi^{k+\frac12}:=\frac12(\psi^{k+1}+\psi^k)$
\begin{thm}[$\ell^2(\bL^2)$ Velocity Estimates]
  \label{thm:velL2rot} Assume that the space dimension is two.  Under
  the assumptions of Theorem~\ref{thm:DIV-est}, the solution $(u_\tau,
  p_\tau)$ of the scheme \eqref{eq:predictor}--\eqref{eq:presupdate}
  in two space dimensions satisfies
\[
\|\ue_\tau -u_\tau\|_{\ell^2(\bL^2)} \leq c \tau^{\frac32}.
\]
\end{thm}
\begin{proof}
  The proof proceeds by a duality argument using the right-inverse
  Stokes operator $S$. By proceeding as in Section~\ref{sub:velL2std}
  we obtain (see \eqref{eq:L2est}):
\[
\|\bar e_\tau \|_{\ell^2(\bL^2)}^2 \le  c(\tau^4 + \|\bar e_\tau-P_\bH\bar e_\tau\|_{\ell^2(\bL^2)}^2).
\]
The estimate \eqref{eq:thm:DIV-est} then immediately implies
\begin{align*}
\|\bar e_\tau \|_{\ell^2(\bL^2)}^2
& \le  c(\tau^4 + \|\DIV\bar e_\tau\|_{\ell^2(\bL^2)}^2) \le c\tau^3.
\end{align*}
Now we observe that
\begin{align*}
  \|e^{k+1}\|_{\bL^2}^2 & = \|\bar e^{k+\frac12} + \tfrac12\delta
  e^{k+1}\|_{\bL^2}^2 \le \tfrac32 \|\bar e^{k+\frac12}\|_{\bL^2}^2 +
  \tfrac34\|\delta e^{k+1}\|_{\bL^2}^2,
\end{align*}
which, along with \eqref{eq:thm:DIV-est}, implies
\begin{align*}
  \|e_\tau\|_{\ell^2(\bL^2)}^2 &  \le \tfrac32 \|\bar e_\tau\|_{\ell^2(\bL^2)}^2 +
  c\|\delta e^{k+1}\|_{\ell^\infty(\bL^2)}^2 \le c\tau^3,
\end{align*}
which completes the argument.
\end{proof}

Let us now show convergence of the velocity in the
$\ell^2(\bH^1)$-norm without any restriction on the space dimension.
\begin{thm}[$\ell^2(\bH^1)$ Velocity Estimates]
\label{thm:velH1rot}
Under the assumptions of Theorem~\ref{thm:DIV-est}, the solution
$(u_\tau, p_\tau)$ of the scheme
\eqref{eq:predictor}--\eqref{eq:presupdate} in two space dimensions
and \eqref{eq:3dpredictor}--\eqref{eq:3dpresupdate} in three space
dimensions satisfies
\[
  \| \bar \ue_\tau  - \bar u_\tau \|_{\ell^\infty(\bH^1)} \leq c \tau.
\]
\end{thm}
\begin{proof}
Observe first that the following holds for all $ k  \in \{0,\ldots,K\}$:
\[
\|\ROT \bar e^{k+\frac12}\|_{\bL^2} \le \sum_{i=1}^k \|\ROT \delta \bar e^{i+\frac12}\|_{\bL^2}
+ \|\ROT \bar e^{\frac12}\|_{\bL^2},
\]
which implies
\[
\|\ROT \bar e^{k+\frac12}\|_{\bL^2} \le  c\tau^{-1}\|\ROT \delta\bar e_\tau\|_{\ell^2(\bL^2)}
+ c \tau  \le c(\tau^{-1} \tau^2 +\tau) \le c \tau,
\]
and, owing to \eqref{eq:grad_eq_divProt}, this concludes the proof
since we have already established that $\|\DIV \bar
e^{k+\frac12}\|_{L^2}\le c \tau^{\frac32}$.
\end{proof}

\begin{rem}[Pressure Error Estimates]
  The same methods and ideas used in Section~\ref{sub:pres} can be
  invoked to show that the pressure satisfies the following estimate
\[
  \| \epsilon_\tau \|_{\ell^2(\Delta)} \leq c \tau.
\]
We omit the details for the sake of brevity,
\end{rem}

\begin{rem}
  Whether Theorem \ref{thm:velL2rot} holds in three space dimensions
  is not clear. The main obstacle in the way is the splitting error
  induced by the splitting of the momentum equation. Based on our
  numerical experiments, we conjecture that both the error estimates in
  Theorem \ref{thm:velL2rot} and Theorem \ref{thm:velH1rot} can be
  improved by a $\tau^{\frac12}$ factor irrespective of the space
  dimension.
\end{rem}

\section{Other Time Marching Techniques}
\label{sec:othertimes}
As mentioned in Remark~\ref{rem:douglas}, the velocity update
\eqref{eq:3dxi}--\eqref{eq:3du} is a sequence of three approximations
of the momentum equation where each approximation consists of
evaluating the second derivative in one of the spatial directions
implicitly with the Crank-Nicolson scheme whereas in the other
directions it either employs the solution from the previous time
level, if no implicit approximation is yet computed in the given
direction, or uses the already computed implicit approximations.  This
observation leads us to propose the following split version of the
second-order backward difference scheme (BDF2) to approximate the
momentum equation:
\[
\frac{3 \eta^{k+1} - 4 u^k +u^{k-1}}{2\tau} - \dxx \eta^{k+1} -
(\dyy+\dzz) u^k + \GRAD p^{\star,k+1} = f^{k+1}.
\]
\[
\frac{ 3\zeta^{k+1} -4 u^k +u^{k-1} }{2\tau} - \dxx \eta^{k+1} - \dyy
\zeta^{k+1} -\dzz u^k+ \GRAD p^{\star,k+1} = f^{k+1}.
\]
\[
\frac{3 u^{k+1} -4 u^k+u^{k-1} }{2\tau} - \dxx \eta^{k+1} - \dyy
\zeta^{k+1} -\dzz u^{k+1} + \GRAD p^{\star,k+1} = f^{k+1}.
\]

We now write the full BDF2 algorithm in a form similar to
\eqref{eq:3dpredictor}--\eqref{eq:3dpresupdate}.  To simplify the
presentation, let us assume that proper approximations of the velocity
and the pressure time derivative are available at $t=-\tau$ and $t=0$.
If these quantities are not available, we start the scheme with a
lower-order approximation at the first time step in order to compute
those approximations.

\begin{itemize}
\item \underline{Pressure predictor:} Denoting by $\pe_0$ the pressure
  field at $t=0$, by $\phi^{\star,0}$ an approximation of
  $ \tau \partial_t \pe(0)$, and by $\phi^{\star,-1}$ an approximation of
  $ \tau \partial_t \pe(-\tau)$ the algorithm is initialized by
  setting $p^{0}=\pe_0$, $\phi^{0}=\phi^{\star,0}$, and $\phi^{-1}=\phi^{\star,-1}$. Then for all
  $k\ge 0$ a pressure predictor is computed as follows:
\begin{equation}
  p^{\star,k+1} =p^{k} + \frac43\phi^{k}-\frac13\phi^{k-1}
\end{equation}
\item \underline{Velocity update:} The velocity update is
  computed by solving the following series of one-dimensional
  problems: Find $\xi^{k+1}$, $\eta^{k+1}$, $\zeta^{k+1}$, and $u^{k+1}$ such that
\begin{align}
  \frac{ 3 \xi^{k+1} - 4 u^k +u^{k-1}}{2 \tau}
  - \LAP u^k + \GRAD p^{\star,k+1} = f^{k+1}, & \quad \xi^{k+1}|_{\partial\Omega} = 0, \\
  \frac{ 3(\eta^{k+1} - \xi^{k+1}) }{2 \tau} - \dxx ( \eta^{k+1} - u^k ) = 0, & \quad \eta^{k+1}|_{x=0,1} = 0,\\
  \frac{3( \zeta^{k+1} - \eta^{k+1} )}{2 \tau} - \dyy ( \zeta^{k+1} - u^k ) = 0, & \quad \zeta^{k+1}|_{y=0,1} = 0,\\
  \frac{ 3(u^{k+1} - \zeta^{k+1} )}{2 \tau} -\dzz ( u^{k+1} - u^k ) =
  0, & \quad u^{k+1}|_{z=0,1} = 0.
\end{align}
\item \underline{Penalty step:} The pressure-correction $\phi^{k+1}$
  is computed by solving
  \begin{equation}
    A \phi^{k+1} = -\frac3{2\tau} \DIV u^{k+1}.
  \end{equation}
\item \underline{Pressure update:} The pressure is updated as follows:
  \begin{equation}
    p^{k+1} = p^{k} + \phi^{k+1} - \chi \DIV \bar u^{k+1}.
  \end{equation}
\end{itemize}
Note that this scheme is formally second-order consistent because
eliminating the intermediate velocities results in a second-order
perturbation of the classical pressure-correction BDF2 scheme.
Numerical experiments show that this algorithm is indeed
unconditionally stable when tested on the unsteady Stokes problem and
its rate of convergence is similar to that of
\eqref{eq:3dpredictor}--\eqref{eq:3dpresupdate}.

\section{Numerical Experiments}
\label{sec:NumExp}
We report in this sections numerical tests aiming at evaluating the
performance of the algorithms
\eqref{eq:predictor}--\eqref{eq:presupdate} in two space dimensions
and \eqref{eq:3dpredictor}--\eqref{eq:3dpresupdate} in three space
dimensions. The space approximation is done using the MAC scheme.

\subsection{Accuracy Tests}
The standard and the rotational versions of the scheme
\eqref{eq:predictor}--\eqref{eq:presupdate} have been tested
numerically on a two dimensional analytic solution and the results
have been reported in \cite{Guermond2010581}. The rate of convergence
with respect to $\tau$ for the velocity in the $\bL^2$-norm for both
versions of the method is about $1.8$ or higher, whereas for the
pressure in the $L^2$-norm it is about $1.85$ for the rotational
version and about $1.5$ for the standard version.

We now investigate the convergence rates in three space dimensions in
$\Omega=(0,1)^3$ using the following solution of the unsteady Stokes
problem (with the appropriate source term):
\begin{align*}
  &\ue_1=(\sin x \cos y \sin z-\sin x \sin y \cos z) \sin t\\
  &\ue_2=(\sin x \sin y \cos z-\cos x \sin y \sin z) \sin t\\
  &\ue_3=(\cos x \sin y \sin z-\sin x \cos y \sin z) \sin t\\
  &\pe=\cos(x+y+z+t).
\end{align*}

\begin{figure}[ht]
  \centerline{
    \includegraphics[width=0.49\textwidth]{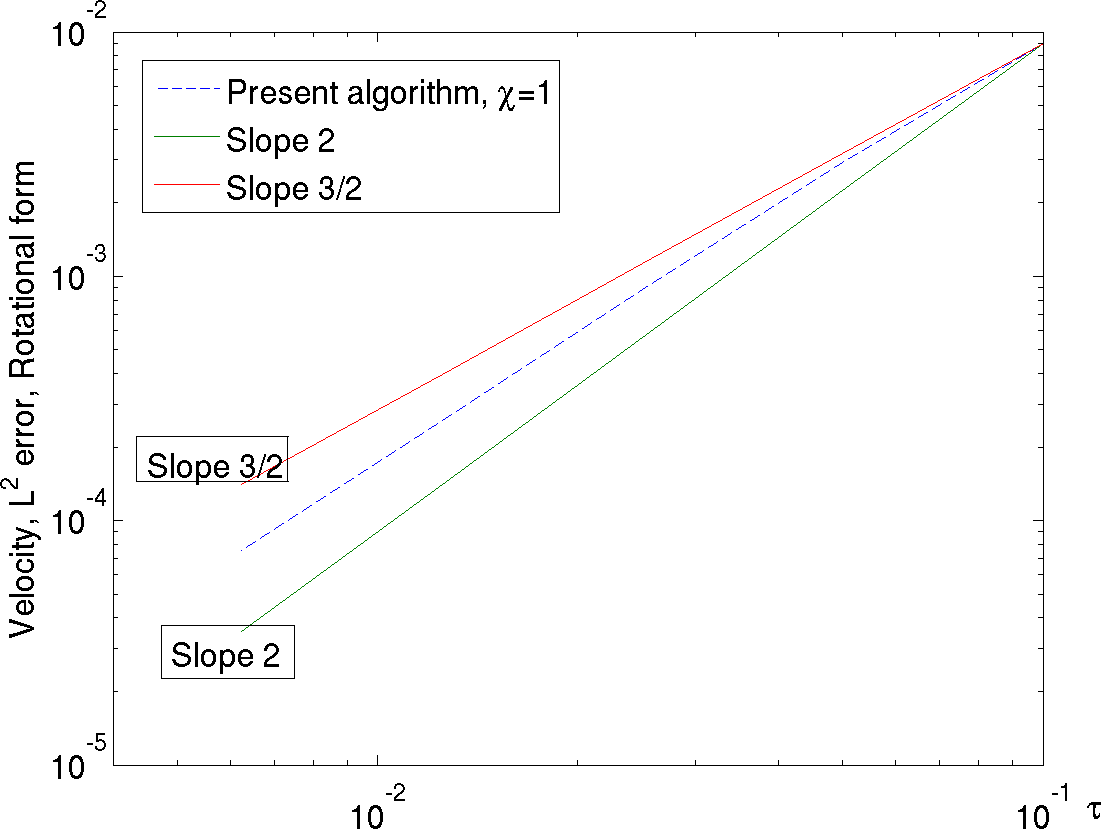}
    \hfill
    \includegraphics[width=0.49\textwidth]{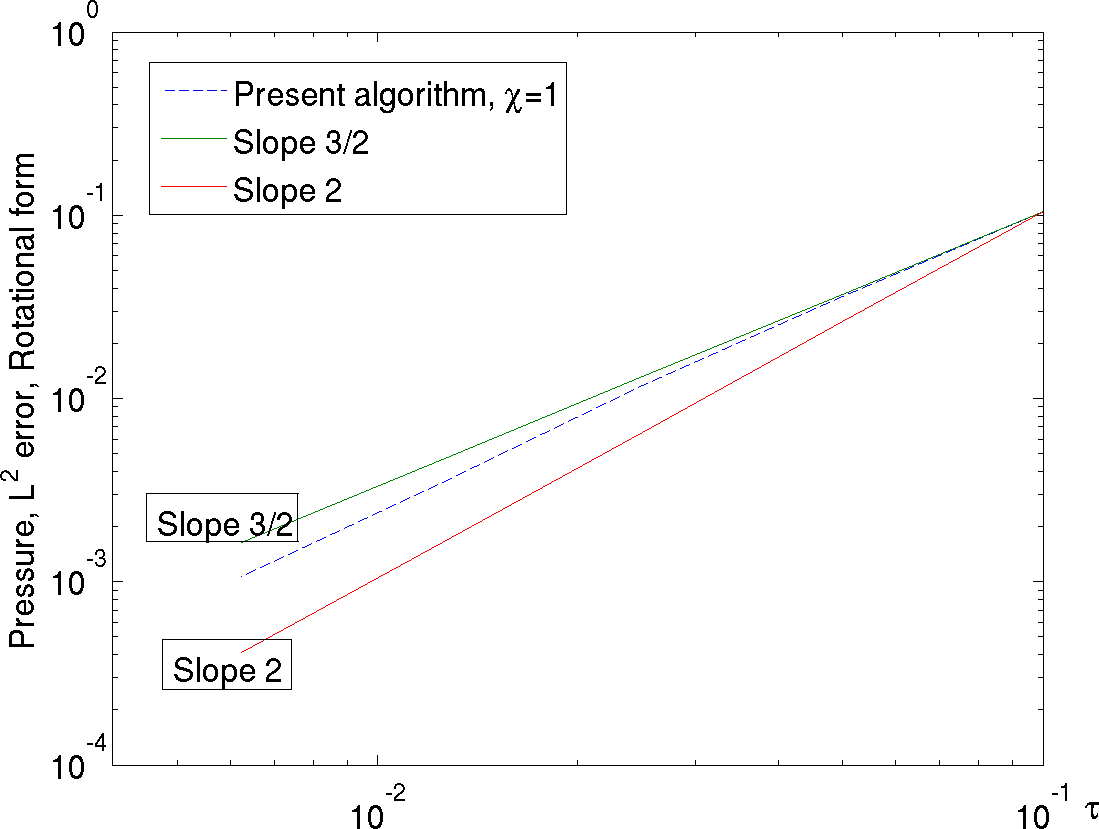}
  }
  \caption{Rotational form ($\chi=1$).  Left: $\bL^2$-norm of
    the error on the velocity (dashed line) at $T=2$ on a uniform grid,
    $100{\times}100 \times 100$; Right: $L^2$-norm of
    the error on the pressure (dashed line).}
\label{errRot}
\end{figure}

We display in the left panel of Figure~\ref{errRot} the $\bL^2$-norm
of the error on the velocity at $T=2$ versus the time step $\tau$ for
the rotational scheme with $\chi=1$.  The $L^2$-norm of the error on
the pressure is displayed in the right panel of the figure.  The
convergence rate on the velocity varies between $1.6$ and $1.8$ while
the convergence rate on the pressure is comprised between $1.5$ and
$1.7$. From tests not reported here, we have observed that the
standard version of the scheme has a convergence rate between $1.6$
and $1.7$ for the velocity and a convergence rate between $1.25$ and
$1.4$ for the pressure. These results suggest that the actual
convergence rates of both schemes are higher than those theoretically
estimated above. However, at the present it is unclear how to improve
these estimates.

\subsection{Splitting vs.~Projection}
\begin{figure}
  \centerline{
    \includegraphics[width=0.49\textwidth]{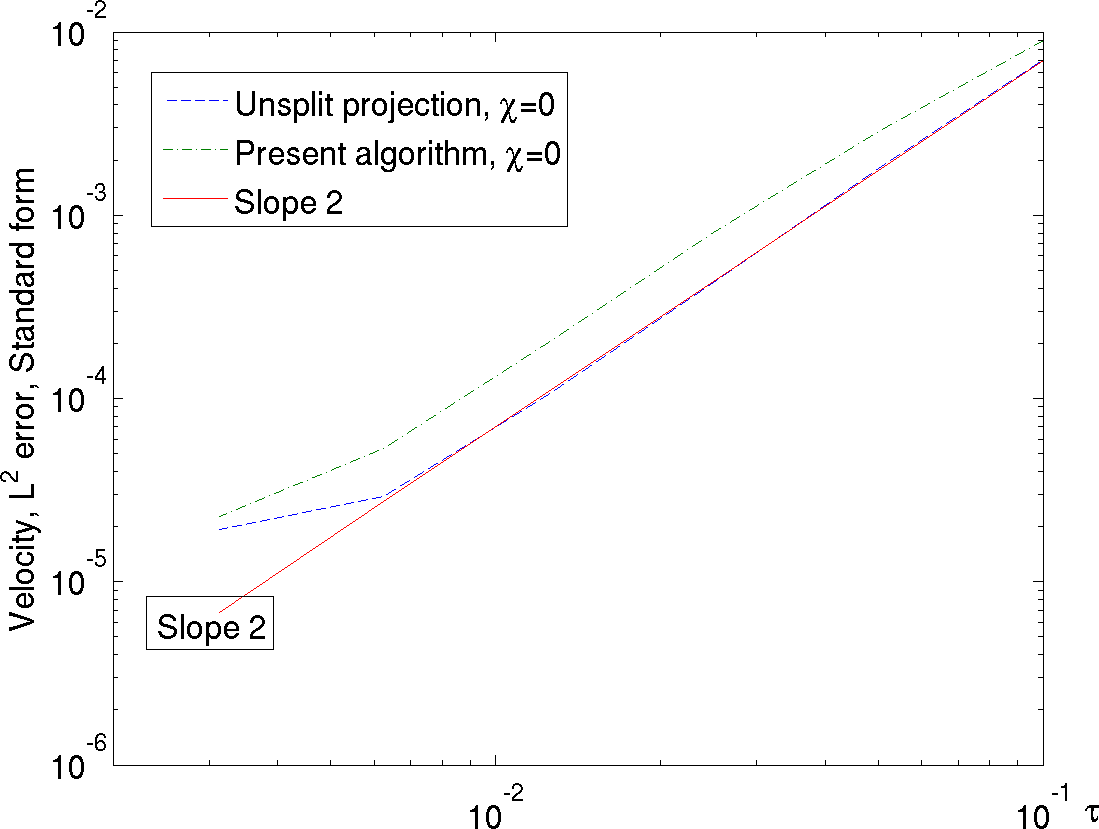}
    \hfill
    \includegraphics[width=0.49\textwidth]{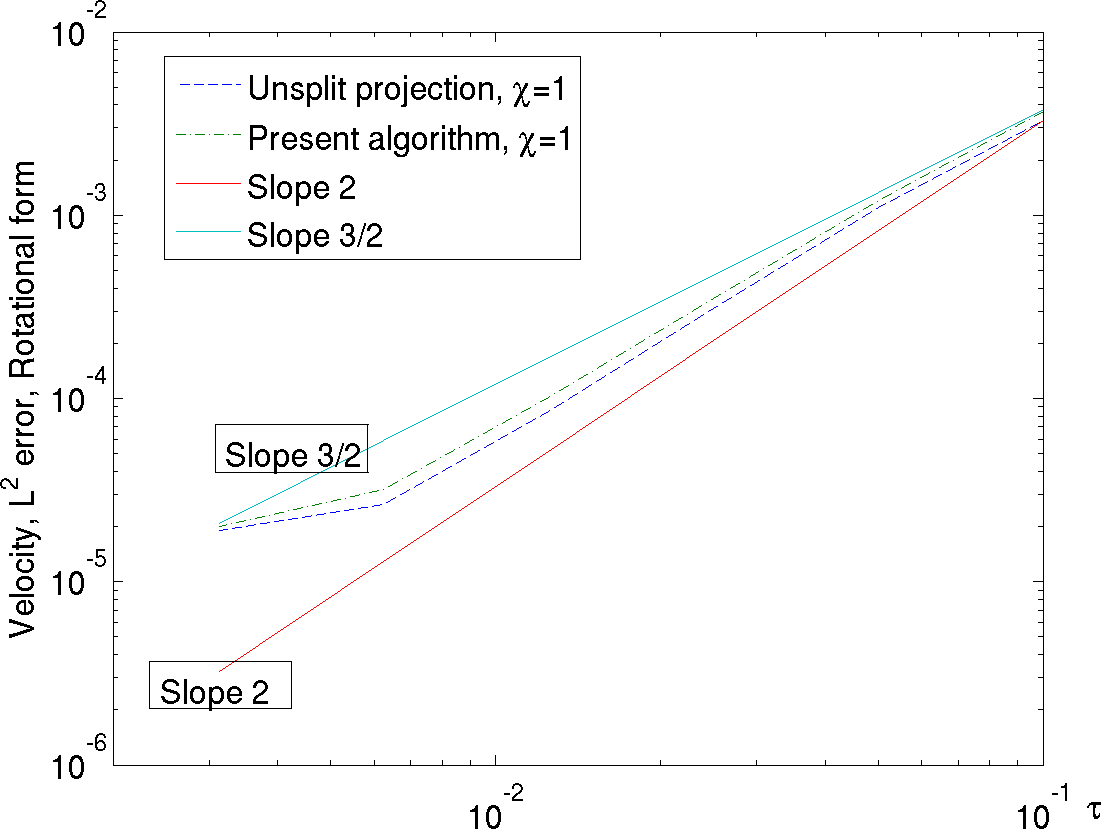}
  }
  \centerline{
    \includegraphics[width=0.49\textwidth]{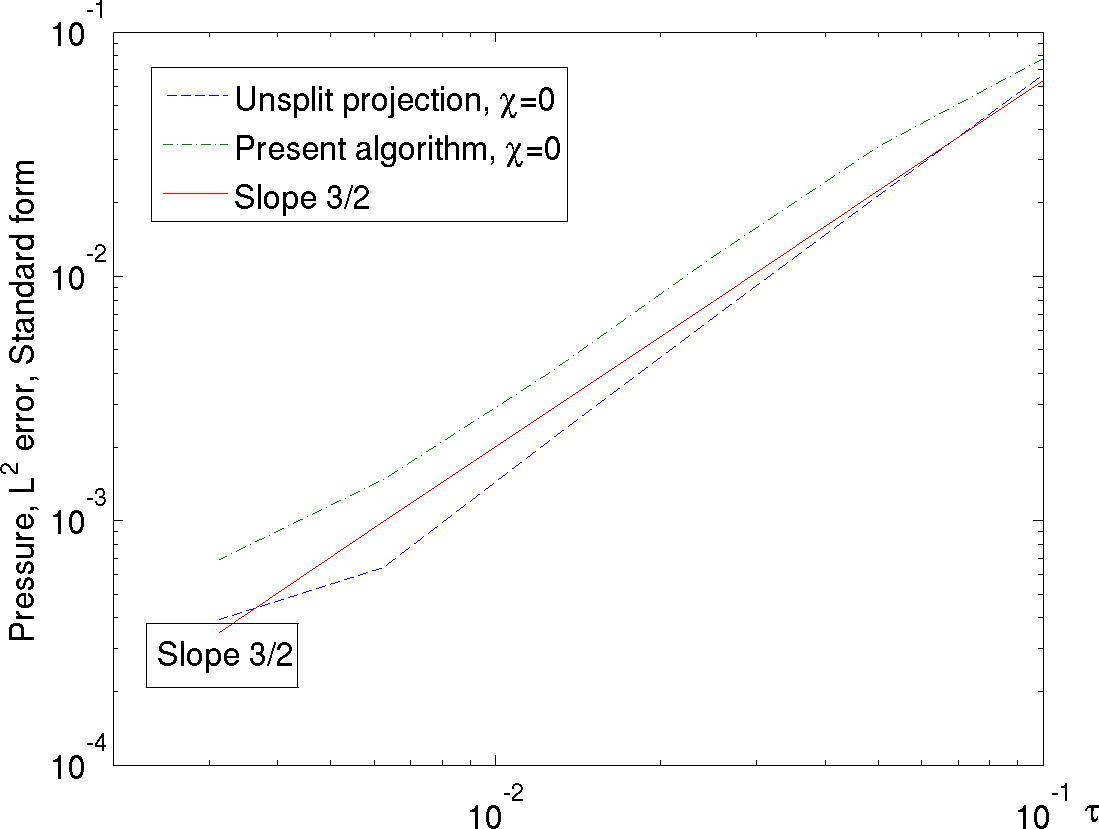}
    \hfill
    \includegraphics[width=0.49\textwidth]{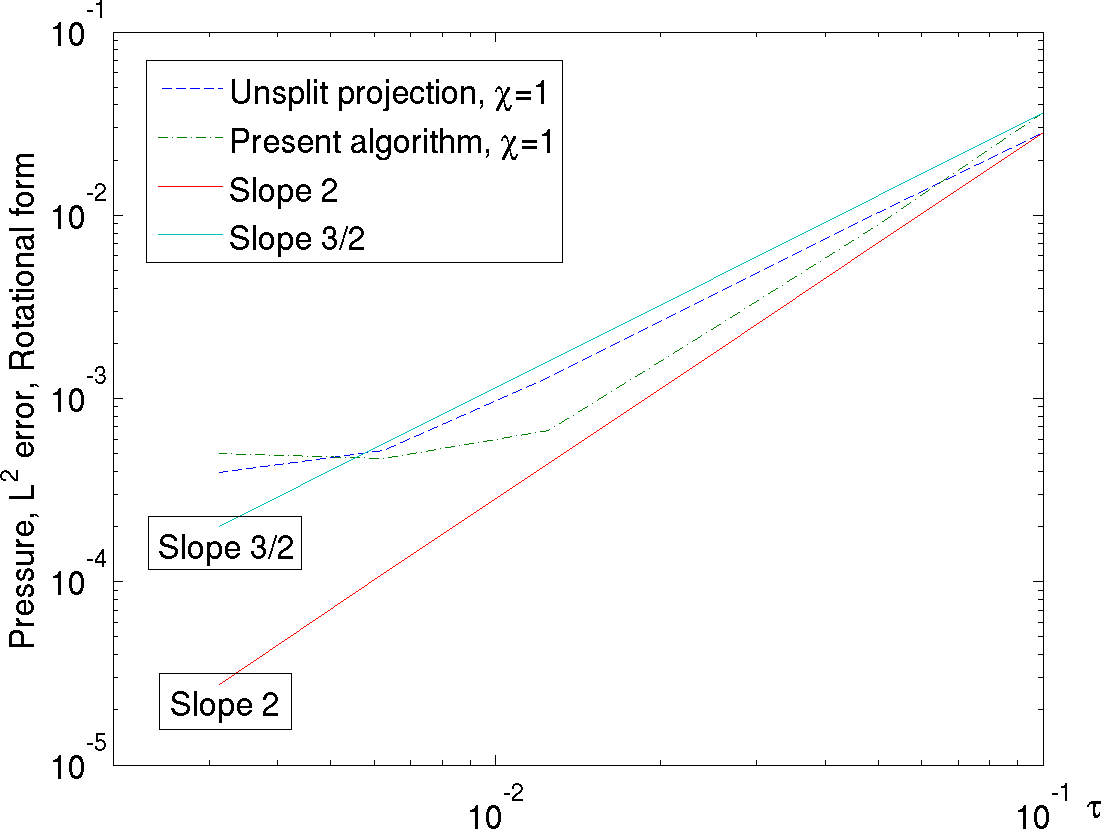}
  }
  \caption{$\bL^2$-norm of the error on the velocity (top) and
    pressure (bottom) at $T=2$ on a uniform $40{\times}40$ grid; Left:
    unsplit projection scheme in standard form (dashed line) and
    scheme \eqref{eq:halfstep}--\eqref{eq:presupdate} with $\chi=0$
    (dash-dotted line).  Right: unsplit projection scheme in
    rotational form with $\chi=1$ (dashed line) and scheme
    \eqref{eq:halfstep}--\eqref{eq:presupdate} with $\chi=1$
    (dash-dotted line).  }
\label{fullVSadi}
\end{figure}

To further illustrate the convergence properties of the present
schemes we now compare it with its unsplit pressure-correction
counterpart, \ie the momentum equation is unsplit and the pressure
correction is computed by solving the Poisson problem ($A=-\LAP_N$).
The comparison is done in two space dimensions in
$\Omega=(0,1)^2$ on the following analytical solution:
\begin{equation}
  \ue= (\sin x \sin(y+t),\cos x \cos(y+t)), \quad \pe=\cos x \sin(y+t).
\label{analytic}
\end{equation}
We show in Figure~\ref{fullVSadi} the error on the velocity and the
pressure as functions $\tau$ for the unsplit second-order projection and the
corresponding results using the present direction splitting schemes.
Clearly, both the standard and the rotational versions of the
direction splitting schemes produce results that are very similar to
those produced by their unsplit counterpart. The largest differences
are observed on the velocity for the standard version of the
schemes. But, even in this case, the direction splitting produces
errors which are only between $1.2$ and $2$ times larger than the errors
produced of the classical standard scheme. The computational
complexity of the present schemes, however, is significantly lower.

\subsection{Lid Driven Cavity}
We compare in this section the performance of the direction splitting algorithm
with its unsplit pressure-correction counterpart on the so-called lid driven cavity.
The computational domain is $\Omega=(0,1)^2$. The boundary conditions
are $u|_{x=0,1,y=0}=0$, $u|_{y=1}=1$ and $v|_{\partial\Omega}=0$.
The computation is done at Reynolds number $R_e=100$
on a MAC grid composed of $40\times40$ nodes and with time step
$\tau=0.01$.  The advection term is computed by means of the
explicit second-order Adams-Bashforth approximation.
The comparison between the two codes is done at $t=1$ and
$t=10$.

\begin{figure}
  \centerline{
    \includegraphics[width=0.49\textwidth]{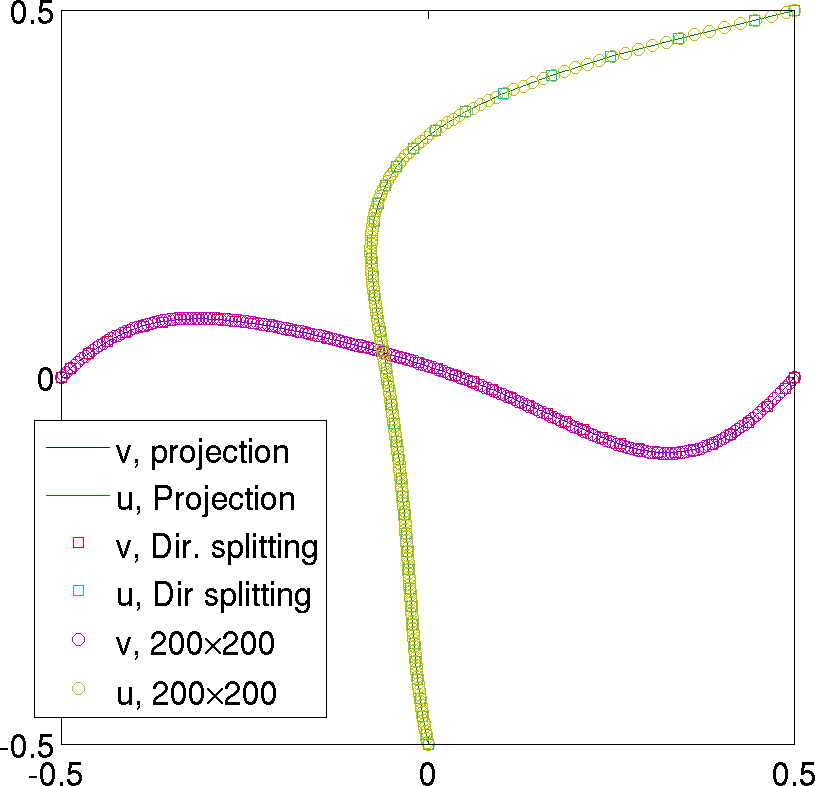}
    \hfill
    \includegraphics[width=0.49\textwidth]{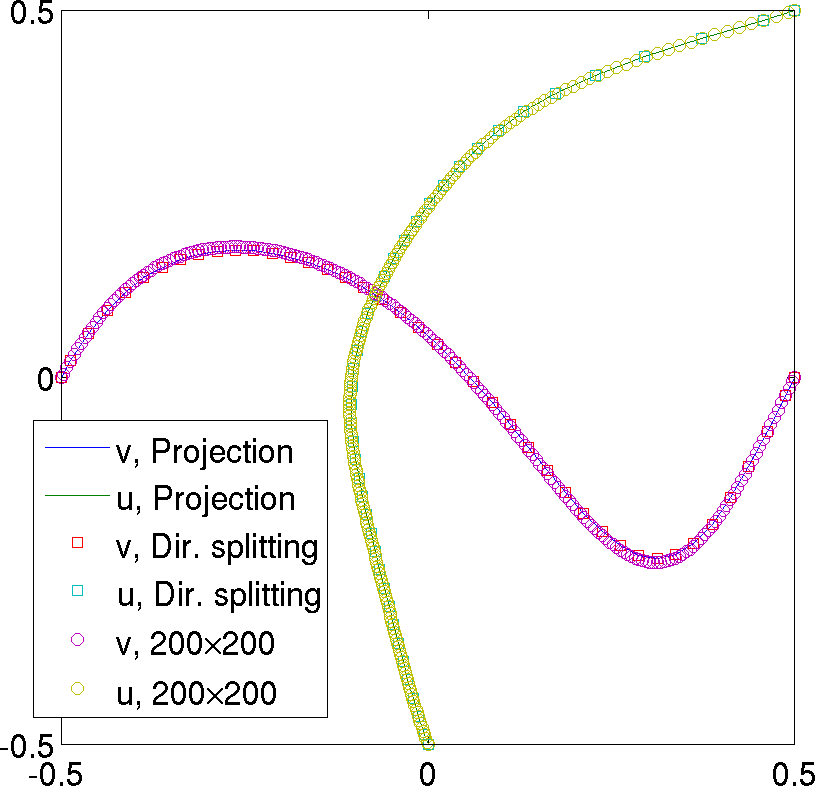}
  }
  \caption{The horizontal (resp. vertical) profiles of the velocity
    alongside the vertical (resp. horizontal) lines through the center of the
    cavity at $R_e=100$ and at $t=1$ (left panel) and $t=10$ (right panel) on
    $40{\times}40$ MAC grid, $\tau=0.01$. Unsplit
    projection scheme (solid line); direction splitting
    scheme ({\tiny$\color{Red}\Box$}, {\tiny$\color{Turquoise}\Box$} symbols); direction splitting on $200{\times}200$ MAC grid with
    $\tau=0.0025$  ($\color{Plum}\circ$, $\color{GreenYellow}\circ$ symbols).}
\label{cavity}
\end{figure}

We show in figure~\ref{cavity} the horizontal and vertical profiles of the
velocity alongside the vertical/horizontal lines through the center of
the cavity.  The results of the two schemes (unsplit and split) are very close to each other;
detailed examination (not reported here) shows that the two sets of results
differ in the fourth decimal digit. For
comparison, we also display with $\circ$ symbols the result of the
split scheme on a MAC grid of $200 \times 200$ nodes and with time step
$\tau=0.0025$.  The three sets of results are visually indistinguishable.

\subsection{Backward Facing Step}
Finally, the new direction splitting method is validated on the
two-dimensional flow over a backward-facing step. Extensive
experimental and computational data on this flow is available in
\cite{ADPS_83} and \cite{Kim_Moin_1985}. Here we compute the solution
to this problem in a rectangular cavity of size $1\times 16$ with a
uniform grid of size $h=0.005$ and a time step $\tau=0.001$. We
prescribe the fully developed parabolic profile with maximum velocity
$\frac32$ at the upper half of the inflow side and we prescribe the
no-slip condition at the lower half. At the outlet we impose
zero-Neumann conditions on the velocity and the zero Dirichlet
condition on the pressure. One important characteristic of the flow is
the length of the recirculation zone behind the step, say $r$. We
report in Table~\ref{Table:1} the results of the present computations
at Reynolds numbers (based on the channel height) $R_e=100$, $200$ and
$400$ and we compare these results with those from
\cite{Kim_Moin_1985}.
% and the experimental findings of \cite{ADPS_83}.
The present scheme yields results
which are in a very good agreement with the existing data.

\begin{table}[h]
  \begin{center}
    \begin{tabular}{||c|c|c||}
    \hline
    $R_e$ & \multicolumn{2}{c||}{$r/s$} \\  \cline{1-3}
            & Current result & Result in \cite{Kim_Moin_1985}\\
             \cline{2-3}
    100 &3.22 &3.2 \\ \hline
    200  &5.33 &5.3 \\ \hline
    400 &8.6 &8.6 \\
    \hline
    \end{tabular}
    \caption{Flow over a backward-facing step. Re-attachment length $r$ divided by the step height $s$ as a function of
            the Reynolds number $R_e$ for the present computations and for the computations of Kim and Moin
            \cite{Kim_Moin_1985}. }
  \end{center}
\label{Table:1}
\end{table}

\subsection{Parallel Implementation}
We have implemented a parallel version of the algorithm
\eqref{eq:3dpredictor}--\eqref{eq:3dpresupdate} with the MAC stencil
using central differences for the first- and second-order
derivatives. The algorithm has been implemented in parallel on a
Cartesian domain decomposition using MPI.
All the one-dimensional
linear systems are solved in parallel using direct solves of the Schur
complement induced by the domain decomposition.  We have verified that
the weak scalability of the code is quasi-perfect up to the maximum
number of processors that were available to us without special request
for allocation, \ie 1024 processors.

Extensive numerical tests have shown that the algorithm is stable
under CFL condition in the Navier-Stokes regime. We have computed a
highly accurate benchmark solution for the start-up flow in a
three-dimensional impulsively started lid-driven cavity of aspect
ratio $1{\times}1{\times}2$ at Reynolds numbers $1000$ and $5000$.
Successive refinements have shown that the velocity field is four digit
accurate at $R_e=5000$ for dimensionless times $t=4$, $8$ and $12$. The
computations have been done in parallel (up to 1024 processors) on
adapted grids of up to 2 billion nodes in three space dimensions. All
these numerical experiments are reported in a forthcoming paper.

\subsection{Further Developments}
We believe that the algorithm presented in the present paper has a lot of
potential for further developments; we are thinking in particular of
academic problems that can be solved in simple geometries with regular grids
\eg simulation of turbulent flows in the atmosphere and in the ocean,
simulation of multiphase flows, stratified flows, variable density
flows, combustion, solution of subgrid problems as part of an
homogenization procedure, \etc

As described in the present paper, the algorithm is suitable only for
simply-shaped domains. However, there are possibilities to impose
boundary conditions either via penalty methods, or fictitious domain
techniques, or via directional adjustment of the grid at the
boundary. The authors have implemented the directional adjustment
procedure and have observed that the resulting scheme is unconditionally
stable and convergent for the time-dependent Stokes problem. These
results will be reported elsewhere.

%============ References ==================
\bibliographystyle{plain}
\bibliography{biblio}

\end{document}